
\magnification=\magstep1

\input amstex

\documentstyle{amsppt}
\overfullrule=0pt
\TagsOnRight

\def\cite#1{{\rm [#1]}}

\leftheadtext{E. MAKAI, Jr. and J. ZEM\'ANEK}

\rightheadtext{ALGEBRAIC ELEMENTS IN A BANACH ALGEBRA}

\topmatter

\title
On the structure of the set of algebraic elements in a Banach algebra 
and their liftings
\endtitle

\author
\leftskip=42.5mm
\baselineskip=10pt
Endre Makai, Jr.\newline
{\eightpoint MTA, Alfr\'ed R\'enyi Institute of Mathematics\newline
H-1364 Budapest, Pf. 127\newline
Hungary\newline
{\lowercase{makai.endre\@renyi.mta.hu\newline
{\rm{http://www.renyi.hu/\~{}makai}}}\newline
\phantom{and}}}\newline
J. Zem\'anek\newline
{\eightpoint Institute of Mathematics\newline
Polish Academy of Sciences\newline
00-656 Warsaw, \'Sniadeckich 8\newline
Poland\newline
\lowercase{zemanek\@impan.pl}}
\endauthor

\leftskip=0pt




\abstract
We generalize earlier results about connected components of idempotents in
Banach algebras, due to B. Sz\H{o}kefalvi Nagy, Y. Kato, S. Maeda,
Z. V. Kovarik, J. Zem\'anek, J. Esterle.
Let $A$ be a unital complex Banach algebra, and $p(\lambda) = \prod\limits_{i
= 1}^n (\lambda - \lambda_i)$ a polynomial over $\Bbb C$, with all roots
distinct.
Let $E_p(A) := \{a \in A \mid p(a) = 0\}$.
Then all connected components of $E_p(A)$ are pathwise connected (locally
pathwise connected) via each of the following three types of paths:
1)~similarity via a finite product of exponential functions (via an
exponential function);
2)~a polynomial path (a cubic polynomial path);
3)~a polygonal path (a polygonal path consisting of $n$ segments).
If $A$ is a $C^*$-algebra, $\lambda_i \in \Bbb R$, let $S_p(A):= \{a\in A \mid
a = a^*$, $p(a) = 0\}$.
Then all connected components of $S_p(A)$ are pathwise connected (locally
pathwise connected), via a path of the form $e^{-ic_mt}\dots e^{-ic_1t}
ae^{ic_1t}\dots e^{ic_mt}$, where $c_i = c_i^*$, and $t \in [0, 1]$ (of the
form $e^{-ict} ae^{ict}$, where $c = c^*$, and $t \in [0,1]$).
For (self-adjoint) idempotents we have by these old papers that the distance
of different connected components of them is at least~$1$.
For $E_p(A)$, $S_p(A)$ we pose the problem if the distance of different
connected components is at least $\min \bigl\{|\lambda_i - \lambda_j| \mid
1 \leq i,j \leq n, \ i \neq j\bigr\}$.
For the case of $S_p(A)$, we give a positive lower bound for these distances,
that depends on $\lambda_1, \dots, \lambda_n$.
We show that several local and global lifting theorems for analytic families
of idempotents, along analytic families of surjective Banach algebra
homomorphisms, from our recent paper with B. Aupetit and M. Mbekhta, have
analogues for elements of $E_p(A)$ and $S_p(A)$.
\endabstract

\endtopmatter

\document
\heading
1. Introduction
\endheading

In this paper all Banach spaces and Banach algebras are over $\Bbb C$, and all
Banach algebras have units whose norms are $1$, and all Banach algebra
homomorphisms preserve the units.
For $X$ a Banach space, $B(X)$ denotes the Banach algebra of all bounded
linear operators $X \to X$.
For a Banach algebra~$A$,
$$
E(A) := \bigl\{a \in A \mid a^2 = a\bigr\},
$$
and for a Banach algebra $A$ with an involution,
$$
S(A) := \bigl\{a \in A \mid a = a^2 = a^* \bigr\}.
$$

We begin with a terminology.
An {\it arc} in a Hausdorff space $X$ is a homeomorphic image of $[0, 1]$ in
$X$.
A {\it path} in a Hausdorff space $X$ is a continuous image of $[0, 1]$.
Since we only consider metric, thus Hausdorff spaces, (locally) arcwise
connectedness is equivalent to (locally) pathwise connectedness.
Namely, Hausdorff continuous images of $[0, 1]$, called {\it Peano continua},
are arcwise connected, cf.\ [Ku], \dots .
Further in this paper we will deal with (locally) pathwise connectedness, and
will look for ``nice'' connecting paths.

The subject began with an observation of B. Sz\H{o}kefalvi Nagy ([SzN42],
Ch.\ \dots \S~3, Hilfssatz, p. 58, [SzN47], Ch.\ \dots, \S~1, 3, p.~350 and
[RSzN] \S~105, Th\'eor\`eme, p.~266) that two orthogonal (i.e., self-adjoint)
projections on a Hilbert space $H$, of distance less than~$1$, are similar via
a unitary.
Later C. Davis [Dav], \dots , Y. Kato [YKa75], Theorem, p.~257, [YKa76],
Theorem~2, p.~367 gave simpler proofs for this theorem.

Y. Kato [YKa76], Ch.\ I, \S~4, 6, p.~33 proved the analogous statement for $X$
a Banach space, for two projections, i.e., idempotent operators in $B(X)$, of
distance less than~$1$, with similarity via an invertible operator.
[YKa76], Ch.\ I, \S~4, 6, Problem 4.13, p.~34 showed that under the same
hypothesis there is an analytic path connecting the two idempotent operators.

S. Maeda [Ma] investigated the set $S(A)$ of self-adjoint projections in
$C^*$-al\-ge\-bras~$A$.
He showed in his Theorem~2 and its Corollary that its connected components are
arcwise connected and locally arcwise connected.
Thus are relatively open among all self-adjoint projections.
He showed in his Lemma~2 that for $e, f \in S(A)$, with $\|f - e\| < 1$, $e$
and $f$ are similar via a self-adjoint involution.
In his corollary to Theorem~2 he showed that two self-adjoint idempotents
belong to the same connected component if and only if they are similar via a
finite product of self-adjoint involutions.
He showed in his Theorem~1 that if $e, f \in S(A)$, are similar via a finite
product of self-adjoint involutions, then they can be connected by a
self-adjoint projection valued path.

Later Z. V. Kovarik [Ko], \S~6, Theorem~1 proved, for $X$ a Banach space, and
$E_0, E_1 \in B(X)$ projections at distance less than $1$ that

\newpage

1)~$E_0, E_1$ can be connected by a projection-valued analytic path of the
form $e^{-itw} E_0 e^{itw}$, $t \in [0, 1]$ and

2)~$E_0, E_1$ can be also connected by a polygonal path consisting of two
segments and

3)~$E_0, E_1$ are similar via an involution.

A consequence of 3) is [Ko], \S~8, Theorem~2: If two projections are connected
by a continuous projection valued path, then they are similar via a finite
product of involutions.

\smallskip
J. Zem\'anek [Ze] investigated the idempotents in Banach algebras.
He obtained that

1)~$e, f \in E(A)$, $\|f - e\| < 1$ implies similarity of $e$ and $f$ via an
exponential, i.e., $f = e^{-c} e e^c$, cf.\ his Lemma 3.1;

2) local arcwise connectedness of $E(A)$ and arcwise connectedness of each
connected component of $E(A)$, cf.\ his Theorem 3.2;

3) $e, f \in E(A)$ are in the same connected component of $E(A)$ if and only
if $f$ is of the form
$$
e^{-c_m} \ldots e^{-c_1} e e^{c_1} \ldots e^{c_m};
$$

4) an idempotent $e$ lies in the centre of $A$ if and only if $\{e\}$ is
isolated in $E(A)$ (i.e., is a connected component of $E(A)$);
if an idempotent $e$ does not lie in the centre of $A$, then the connected
component of $E(A)$ containing $e$ contains a (complex) line, hence is
unbounded.

J. Esterle [Es], Theorem and its proof proved refinements of J. Zem\'anek's
results.
He obtained that

1) $e, f \in E(A)$, $\|f - e\| < 1$ implies similarity of $e$ and $f$ or the
form
$$
f = e^{-c''} e^{-c'} e e^{c'} e^{c''} = (1 - c'')(1 - c')e
e(1 + c') (1 + c''), \text{ where } (c')^2 = (c'')^2 = 0,
$$
where even one of the factors containing $c'$ can be omitted, leaving this
formula valid;

2) pathwise connectedness of each connected component of $E(A)$, by polynomial
paths (and at the same time similarities via a finite product of exponential
functions) of the form
$$
\aligned
&e^{-c_n'' t} e^{-c_n't} \ldots e^{-c_1''t} e^{-c_1't} e e^{c_1't} e^{c_1''
t} \ldots  e^{c_n't} e^{c_n'' t}\\
&=(1 - c_n''t)(1 - c_n't) \ldots (1 - c_1'' t)(1 - c_1't) e(1 + c_1't)(1 +
c_1'' t)\ldots (1 + c_n't)(1 + c_n''t),
\endaligned
$$
where $(c_i')^2 = (c_i'')^2 = 0$, for $1 \leq i \leq n$, where one of the
factors containing $c_1'$ can be omitted leaving this formula valid.
Here, for $\| f - e\| < 1$ we have $n = 1$, thus 

\newpage

we have a connection via a
cubic polynomial path.

\smallskip
For some further results about $E(A)$ cf.\ the references in this paper.

Now we turn to the other subject of our paper.
For unital complex Banach algebras $A$ we write
$$
E(A) := \{a \in A \mid a^2 = a\}.
$$
If $A$ is a unital complex Banach algebra with continuous involution ${}^*$,
then we write
$$
S(A) := \{a \in A \mid a^2 = a = a^*\}.
$$
Let $\pi: B \to A$ be a unit preserving homomorphism between two unital
complex Banach algebras $B$ and $A$.
If $A$ and $B$ are Banach algebras with continuous involutions ${}^*$, then we
additionally suppose that $\pi$ is involution-preserving.
Then clearly
$$
\pi E(B) \subset E(A),
$$
and for unital complex Banach algebras with continuous involutions
$$
\pi S(B) \subset S(A).
$$

From now on suppose that $\pi$ is surjective.
We say that the {\it lifting property holds for} $\pi: B \to A$, if
$$
\pi E(B) = E(A), \ \text{ or } \ \pi S(B) = S(A), \ \text{ respectively.}
$$

We write for Banach spaces $X$, $Y$ \ $\Cal B(X, Y)$ for the Banach space of
all bounded linear operators $X \to Y$.
We write $\Cal B(X) := \Cal B(X,X)$, and we write $\Cal K(X)$ for the Banach
space of compact linear operators in $\Cal B(X)$.

For $H$ a Hilbert space, and $\pi : \Cal B(H) \to \Cal B(H) / \Cal K(H)$ the
canonical mapping, we have $\pi S(\Cal B(H)) = S(\Cal B(H) / \Cal K(H))$, cf.\
[Ca], Theorem 2.4 and [dlH], Proposition~7.
For any Banach algebra $A$ and $\pi$ the canonical mapping $A \to A/\text{\rm
rad }A$ ($\text{\rm rad}(\cdot)$ being the radical)
we have $\pi E(A) = E(A / \text{\rm rad }A)$, cf.\ [Ri], Theorem 2.3.9 and
[IKa], p.~125.
An analogue of the above mentioned first result is $\pi E(\Cal B(H)) = E(\Cal
B(H) / \Cal K(H))$, which is due to [La]; we are grateful to Prof.\ J.-Ph.\
Labrousse for personally explaining the difficult passages of his paper.
In fact [La] proved more.
Suppose $U \subset \Bbb C$ with $0 \in U$ is open, and let us have an analytic
map $q: U \to E(\Cal B(H) / \Cal K(H))$.
Then there exist $V \subset \Bbb C$ open, such that $0 \in V \subset U$, and
an analytic map $p: V \to E(\Cal B(H))$, such that $\pi(p(\lambda)) =
q(\lambda)$ for each $\lambda \in V$.
This is called a {\it local lifting theorem}.
If we can choose $V = U$, then we speak about a {\it global lifting theorem}.

In [AMMZ03] there are proved several further local and global lifting
theorems, under hypotheses that the spectra of all elements of $\text{\rm Ker
}\pi$ are ``small''.
Observe 

\newpage

that the spectra of compact operators on a Banach space $X$ are either
finite or are of the form $\{0\} \cup \{\lambda_n \mid n \in \Bbb N\}$, where
$\lambda_n \to 0$, thus they are ``small''.
So no compactness of operators is necessary, but only ``small spectra''.

All these results of [AMMZ03] were strengthened in [AMMZ14].
There not only the idempotents were analytic functions of $\lambda \in U$ (for
$u \subset \Bbb C$ open), but also the surjective unit-preserving Banach
algebra homomorphisms, written as $\pi(\lambda)$.
Under the strongest spectral hypothesis that the spectra of all elements of
all kernels $\text{\rm Ker }\pi(\lambda)$ are $\{0\}$, we proved global
lifting theorems, both for $E(\cdot)$ and for $S(\cdot)$.
For the case of $S(\cdot)$, both the idempotents $q(\lambda)$ and the
${}^*$-homomorphisms $\pi(\lambda)$ were real analytic maps from some open
subset $G$ of $\Bbb R$ with $0 \in G$.
Actually in [AMMZ14] not only a single analytic family of idempotents could be
lifted, but even a mutually orthogonal sequence of analytic families of
idempotents could be lifted to a mutually orthogonal sequence of analytic
families of idempotents.
Here two idempotents $e, f$ in some Banach algebra are {\it orthogonal}, if
$ef = fe = 0$.

Under weaker spectral hypotheses, [AMMZ14] proved local lifting theorems, both
for $E(\cdot)$ and $S(\cdot)$.
For the unital complex Banach algebra case it was sufficient to suppose that
the spectra of all elements of $\text{\rm Ker }\pi(0)$ did not
disconnect~$\Bbb C$.
For the case of unital complex Banach algebras with continuous involutions
(with real analycity of $q(\cdot)$ and the ${}^*$-homomorphisms $\pi(\cdot)$
like above) it was sufficient to suppose that the spectra of all elements of
$\text{\rm Ker }\pi(0)$ was totally disconnected (i.e., they did not contain
any connected subsets consisting of more than one points; this property
implies that they did not disconnect~$\Bbb C$).
Observe that in both of these cases we had no hypotheses for the spectra of
elements of $\text{\rm Ker }\pi(\lambda)$, for $\lambda \neq 0$.

\smallskip
A large part of the results of this paper have been announced in [MZ].


\heading
2. Theorems
\endheading


Let $F$ be a commutative field, and $A$ a unital $F$-algebra.
In general we will write $0$ or $1$ both for the zero or unit in $F$ and $A$,
but it will be clear which is meant; but sometimes we will write $0_F$ and
$0_A$, and $1_F$ and $1_A$.
For $n \geq 2$ an integer, and $\lambda_1, \dots, \lambda_n \in F$ distinct we
write
$$
p(\lambda) := \prod\limits_{i = 1}^n (\lambda - \lambda_i),
$$
that is a polynomial over $F$, and
$$
\aligned
E(A) &:= \{a \in A \mid a = a^2\} \ \text{ and }\\
E_p(A) &:= \{a \in A \mid p(a) = 0_A\}.
\endaligned
$$

\newpage

We consider $n$ and $\lambda_1, \dots, \lambda_n$ as fixed, in the whole paper.
A particular case is when $A$ is the algebra of all bounded linear operators
on a complex vector space $X$ (with $F = \Bbb C$).

We say that
$$
\{e_1, \dots, e_n\} \subset E(A) \ \text{ forms a {\it partition of unity}}.
$$
If $e_i e_j = 0$ for $1 \leq i,j \leq n$, $i \neq j$, and $\sum\limits_{i =
1}^n e_i = 1$.
In other words, we have an $E(A)$-valued measure on $F$, concentrated on
$\{\lambda_1, \dots, \lambda_n\} \subset F$, with total mass $1_A$, with the
measure of $\{\lambda_i\}$ being~$e_i$.

\proclaim{Proposition 1}
With the above notations, we have for $a \in A$ that
$$
a \in E_p(A)
$$
if and only if $a$ is of the form
$$
a = \sum_{i = 1}^n \lambda_i e_i, \text{ where } \{e_1, \dots, e_n\} \subset
E(A) \text{ is a partition of unity}.
$$
Here, for each $1 \leq i \leq n$, the function $a \mapsto e_i = e_i(a)$ is a
polynomial, with coefficients in~$F$.
These polynomials only depend on $\lambda_1, \dots, \lambda_n$.

If $a = \sum\limits_{i = 1}^n \lambda_i e_i = \sum\limits_{j = 1}^m \mu_j
f_j$, for $\{e_1, \dots, e_n\}, \{f_1, \dots, f_m\} \subset E(A)$ partitions
of unity, with distinct $\lambda_i$'s and distinct $\mu_j$'s, and with
$e_i \neq 0$, $f_j \neq 0$ for $1 \leq i \leq n$ and $1 \leq j \leq m$, then
$m = n$, and, after a permutation of the indices, $\lambda_i = \mu_i$ and $e_i
= f_i$, for $1 \leq i \leq n$.
In particular, for $\{e_1, \dots, e_n\}$, $\{f_1, \dots, f_n\} \subset E(A)$
partitions of unity, with distinct $\lambda_i$'s, $a = \sum\limits_{i =
1}^n \lambda_i e_i = \sum\limits_{i = 1}^n \lambda_i f_i$ implies $e_i = f_i$
for $1 \leq i \leq n$.
\endproclaim

\proclaim{Corollary 2}
Let $A = B(X)$, where $X$ is a complex Banach space.
Then for $T \in B(X)$ we have $T \in E_p(A)$ if and only if there exists a
direct sum decomposition $X_1 \oplus \ldots \oplus X_n$ of $X$, where each
$X_i$ is a closed subspace of $X$, such that
$$
T \mid X_i = \lambda_i \cdot id_{X_i}, \ \text{ for each } \ 1 \leq i \leq n.
$$
\hfill $\blacksquare$
\endproclaim

\newpage

{\bf{Remark A.}}
Now we show by an example that our considerations do not work if some
$\lambda_i$'s are equal, even in the ``simplest'' case of $A = B(H)$, the
algebra of bounded linear operators of a Hilbert space.
We consider the simplest such polynomial: $p(\lambda) := \lambda^2$.
Then we can write operators $T \in B(H \oplus H)$ in $2 \times 2$ block matrix
form.
If $T$ is superdiagonal, i.e.,
$$
T = \left( \matrix 0 & T_{12}\\ 0 & 0\endmatrix\right),
$$
we have $T^2 = 0$, and the structure of $T$ can be as complicated as the
structure of $T_{12} \in B(H)$ can be.

\vskip.1cm

{\it If moreover $F$ and $A$ have involutions ${}^*$, we always suppose}
$\lambda_i^* = \lambda_i$, $1 \leq i \leq n$, and we write
$$
\aligned
S(A) &:= \{a \in A \mid a = a^2 = a^*\}, \ \text{ and }\\
S_p(A) &:= \{a \in A \mid a = a^* , \ p(a) = 0_A\}.
\endaligned
$$

For the $C^*$-algebra case, $\lambda_i \in \Bbb R$ is a natural hypothesis.
In fact, if $q(\lambda) := \pi \{(\lambda - \lambda_i) \mid 1 \leq i \leq
n, \lambda_i \in \Bbb R\}$, then $p(a) = 0$ and $a = a^*$ imply $q(a) = 0$.
Thus we could use $q(\cdot)$ rather than $p(\cdot)$.

A particular case is when $A$ is the algebra of all bounded linear operators
on a complex Hilbert space $H$ (with $F = \Bbb C$).

Then we have the analogue of Proposition~1.

\proclaim{Proposition 3}
With the above notations, we have for $a \in A$ that
$$
a \in S_p(A)
$$
if and only if $a$ is of the form
$$
a = \sum_{i = 1}^n \lambda_i e_i,  \text{ where }  \{e_1, \dots, e_n\} \subset
S(A) \text{ is a partition of unity}.
$$
Here, for each $1 \leq i \leq n$, the function $a \mapsto e_i = e_i(a)$ is a
polynomial, with coefficients in $F$, all of which are self-adjoint.
These polynomials only depend on $\lambda_1, \dots, \lambda_n$, and coincide
with those from Proposition~1 (except that here $\lambda_i^* = \lambda_i$).

Of course, the last statement of Proposition~1 remains valid if we replace
$E(A)$ by $S(A)$, and require $\lambda_i^* = \lambda_i$, $\mu_j^* = \mu_j$,
for $1 \leq i \leq n$, $1 \leq j \leq m$.
\endproclaim

\proclaim{Corollary 4}
Let $A = B(H)$, where $H$ is a complex Hilbert space.
Then for $T \in B(H)$ we have $T \in S_p(A)$ if and only if there exists an
orthogonal direct sum decomposition $H_1 \oplus \ldots \oplus H_n$ of $H$,
where each $H_i$ is a closed subspace of $H$, such that
$$
T \mid H_i = \lambda_i \cdot id_{H_i}, \ \text{ for each } \ 1 \leq i \leq n.
$$
\hfill $\blacksquare$
\endproclaim

\newpage

\proclaim{Corollary 5}
Let $F$ be a commutative topological field, and $A$ a unital topological
$F$-algebra.
Then the function $a\mapsto e_i(a)$ from Propositions 1 and 3, being a
polynomial, is continuous, for each $1 \leq i \leq n$.\hfill $\blacksquare$
\endproclaim

{\bf{Remark B.}}
For an ordered partition of unity $(e_1, ..., e_n)$, the function $(e_1, ...,
e_n) \mapsto \sum\limits_{i = 1}^n \lambda_i e_i \in E_p(A)$ is clearly
continuous.
By Corollary~5, for $a \in E_p(A)$, the function $a \mapsto (e_1(a), \dots,
e_n(a))$ is continuous.
Their composition $a \mapsto (e_1(a), ..., e_n(a)) \mapsto \sum \lambda_i
e_i(a)$ is identity by Proposition~1.\
We show that the other composition\break
$(e_1, \dots , e_n) \mapsto \sum \lambda_i
e_i \mapsto \Bigl(e_1\bigl(\sum \lambda_i e_i\bigr), \dots,
e_n \bigl(\sum\lambda_i e_i\bigr)\Bigr)$ is also identity.
We calculate, e.g.,
$e_1\bigl(\sum \lambda_i e_i\bigr)$.
Analogously, as in the proof of Proposition~1, {\bf 2}, we have
$$
e_1\Bigl(\sum \lambda_i e_i\Bigr) = \prod_{j = 2}^n \biggl( \sum_{i = 1}^n
(\lambda_i - \lambda_j)e_i\biggr) \biggm/ \prod_{j = 2}^n (\lambda_1
- \lambda_j)
$$
like in that proof, this product equals
$$
e_1^{n - 1} \prod_{j = 2}^n (\lambda_i - \lambda_j) \Bigm/ \prod_{j = 2}^n
(\lambda_1 - \lambda_j) + \sum_{i = 2}^n e_i^{n - 1} \prod_{j = 2}^n
(\lambda_i - \lambda_j) \Bigm/ \prod_{j = 2}^n (\lambda_1 - \lambda_j) = e_1.
$$
Hence these maps constitute two homeomorphisms, inverse to each other.

\vskip.1cm

The analogous statement holds for ordered self-adjoint partitions of unity,
and $a \in S_p(A)$, by Proposition~3, and the above arguments.

These statements mean that when investigating local connectedness or
connectedness of components of $E_p(A)$ or $S_p(A)$ via analytic or polynomial
paths, we obtain the same answers for connected components of ordered
(self-adjoint) partitions of unity.
For this observe only that the maps $(e_1, \dots, e_n) \mapsto \sum\limits_{i
= 1}^n \lambda_i e_i$ and $a \mapsto (e_1(a), \dots, e_n(a))$ are polynomial
(hence analytic) maps.
Of course, for polygonal maps we do not have an equivalence, but as we will
see in the prooofs of Theorems 12, 13, 14, we will deduce polygonal
connections in $E_p(A)$ and $S_p(A)$ from polygonal connections in ordered
partitions of unity, and ordered self-adjoint partitions of unity.
For this observe that the map $(e_1, \dots, e_n)\mapsto \sum\limits_{i =
1}^n \lambda_i e_i$ is linear.

Later in this paper we will only consider $E_p(A)$ and $S_p(A)$.

Apart from Theorems \dots and \dots, from now on, in the whole paper we
restrict 

\newpage

our attention to the case when $A$ is a unital complex Banach
algebra, or sometimes a $C^*$-algebra, or more generally a unital complex
Banach algebra with a continuous involution.
For $x \in A$ we write $\sigma(x)$ for the {\it spectrum of} $x$ in~$A$.

{\bf{Remark C.}}
For unital complex Banach algebras $A$ Corollary~5 can be shown also in other
ways.

1) For $a \in E_p(A)$, by the spectral mapping theorem, we have
$$
\emptyset \neq \sigma (a) \subset \{\lambda_1, \dots, \lambda_n\}.
$$
Then $e_i$, for $1 \leq i \leq n$, can be obtained as a Riesz idempotent
$$
e_i := \frac{1}{2\pi i} \int\limits_{\Gamma_i} (a - \lambda)^{-1} d\lambda,
\tag{C.1}
$$
where $\Gamma_i$ is a small circle with centre $\lambda_i$, and for $i \neq j$
we have $\Gamma_i \cap \Gamma_j = \emptyset$.
Then also $a = \sum\limits_{i = 1}^n \lambda_i e_i$.
Now (C.1) also implies continuity of the function $a \mapsto e_i = e_i(a)$
(and self-adjointness of $e_i$ for $C^*$-algebras, and $\lambda_i$'s real).

2) Let $a$, $a + \varepsilon \in E_p(A)$, where $\|\varepsilon \| \leq 1$, and
$1 \leq i \leq n$.
Then
$$
\align
\| e_i (a + \varepsilon) - e_i(a)\| &= \|p_i(a + \varepsilon) - p_i(a)\| \\
&= \left\| \left[ \prod_{j = 1 \atop j \neq i}^n (a + \varepsilon - \lambda_j)
- \prod_{j = 1 \atop j \neq i}^n (a - \lambda_j)\right] \Bigm/ \prod_{j =
1 \atop j \neq i}^n (\lambda_i - \lambda_j)\right\|\\
&\leq \left[ \prod_{j = 1 \atop j \neq i}^n (\|a\| + \|\varepsilon\| +
|\lambda_j|) - \prod_{j = 1\atop j \neq i}^n (\|a\| +
|\lambda_j|)\right] \Bigm/ \prod_{j = 1\atop j\neq i}^n |\lambda_i
- \lambda_j|\\
&= \sigma_{\|a\|} (\|\varepsilon\|),
\endalign
$$
where $\sigma_{\|a\|}$ means that the constant in the $\sigma$ sign depends on
$\|a\|$ (observe that $\lambda_1,\dots, \lambda_n$ are fixed in the whole
paper).

\vskip.1cm

{\bf{Remark D.}}
For $A$ a $C^*$-algebra (which is the only case of ${}^*$-algebras that we
will be able to handle when investigating paths) it is no restriction of
generality that we restricted our attention to real $\lambda_i$'s.
In fact, let $a = a^* \in A$.
Then $\sigma(a) \subset \Bbb R$, cf.\ [En], p.~1.
Therefore for $\lambda \in \Bbb C \setminus \Bbb R$ \ $a - \lambda$ is
invertible.
So if we had admitted also non-real $\lambda_i$'s, we could have omitted all
the factors $a - \lambda_i$ with $\lambda_i \in \Bbb C \setminus \Bbb R$  from

\newpage

$p(\lambda)$, thus obtaining a polynomial $q(\lambda)$ with $q(a) = 0$.

\vskip.1cm

\proclaim{Theorem 6}
Let $A$ be a unital complex Banach algebra.
Then $E_p(A)$ is locally polynomially connected, via cubic polynomial paths,
for $t \in [0,1]$.

Moreover, suppose that $a_0, a_1 \in E_p(A)$, $a_0$ is fixed, and $\|a_1 -
a_0\|$ is sufficiently small.
Then $a_0$ and $a_1$ are similar, via some exponential (thus invertible
element), i.e.,
$$
a_1 = e^{-c} a_0 e^c.
$$
This implies an analytic path $a(t) \in E_p(A)$ for $t \in [0, 1]$, from $a_0$
to $a_1$, namely
$$
a(t) = e^{-ct} a_0 e^{ct}.
$$
The distance of this path from $a_0$ tends to $0$, if $\|a_1 - a_0\| \to 0$.

Hence the connected component of $a_0$ in $E_p(A)$ is locally pathwise
connected, via similarity with an exponential function.
\endproclaim

\proclaim{Theorem 7}
Let $A$ be a unital complex Banach algebra, and let $c$ be a connected
component of $E_p(A)$.
Then $C$ is a relatively open subset of $E_p(A)$.
Let $a_0, a_1 \in C$.
Then $a_0$ and $a_1$ are similar, via a finite product of exponentials (that
is invertible), i.e.,
$$
a_1 = e^{-cm} \ldots e^{-c_1} a_0 e^{c_1} \ldots e^{c_m},
$$
for some integer $m \geq 1$, where $c_i \in A$.
This implies an analytic path $a(t) \in C$ from $a_0$ to $a_1$, for $t \in [0,
1]$, namely
$$
a(t) = e^{-c_m t} \ldots e^{-c_1 t} a_0 e^{c_1 t} \ldots e^{c_m t}.
$$
Additionally, we may suppose
$$
c_1^2 = \dots = c_m^2 = 0,
$$
which implies a polynomial path
$$
\tilde a(t) = e^{-c_m t} \ldots e^{-c_1 t} a_0 e^{c_1 t} \ldots e^{c_m t} = (1
- c_m t) \ldots (1 - c_1 t) a_0(1 + c_1 t) \ldots (1 + c_m t),
$$
from $a_0$ to $a_1$, for $t \in [0, 1]$.
Here one of the factors containing $c_1$ can be deleted, without changing the
value of the right-hand side in the last equation.

Hence $c$ is pathwise connected via similarities with finite products of
exponential  functions, and also with polynomial paths.
Moreover, there is a single path satisfying both properties.
\endproclaim

\newpage

\proclaim{Corollary 8}
Let $A$ be a unital complex Banach algebra.
Let $a_0 \in E_p(A)$.
Then $a_0$ is in the centre of $A$ if and only if the connected component of
$a_0$ in $E_p(A)$ is $\{a_0\}$ (i.e., $a_0$ is isolated in $E_p(A)$).
\endproclaim

\proclaim{Theorem 9}
Let $A$ be a unital complex Banach algebra.
If some connected component of $E_p(A)$ does not intersect the centre of $A$,
then any element of $C$ is contained in a complex line entirely contained
in~$C$.
In particular, $C$ is unbounded.
\endproclaim

\proclaim{Corollary 10}
Let $A$ be a unital complex Banach algebra.
Then $E_p(A)$ is a union of its isolated point and of complex lines.
\endproclaim

\proclaim{Theorem 11}
Let $\lambda_1, \dots, \lambda_n \in \Bbb R$.
Let $A$ be a unital complex $C^*$-algebra.
Suppose that $a_0, a_1 \in S_p(A)$, $a_0$ is fixed, and $\|a_1 - a_0\|$ is
sufficiently small.
Then $a_0$ and $a_1$ are similar via a unitary,
$$
a_1 = e^{-ic} a_0 e^{ic}, \ \text{ where } \ c = c^* \in A, \ \text{ and
$\|c\|$ is small}.
$$
This implies a self-adjoint analytic path $a(t) \in S_p(A)$ from $a_0$ to
$a_1$, for $t \in [0, 1]$, namely
$$
a(t) = a(t)^* = e^{-ict} a_0 e^{ict}.
$$
The distance of this path from $a_0$ tends to $0$, if $\|a_1 - a_0\| \to 0$.

Hence the connected component of $a_0$ in $S_p(A)$ is locally pathwise
connected via a similarity with a unitary valued exponential function.
\endproclaim

\proclaim{Theorem 12}
Let $\lambda_1, \dots, \lambda_n \in \Bbb R$.
Let $A$ be a unital complex $C^*$-algebra, and let $C$ be a connected
component of $S_p(A)$.
Then $C$ is a relatively open subset of $S_p(A)$.
Let $a_0, a_1 \in C$.
Then $a_c$ and $a_1$ are similar.
Via some unitary $s \in A$, i.e.,
$$
a_1 = s^{-1} a_0 s = s^* a_0 s_s \ \text{ where } \ s^{-1} = s^*.
$$
Here, for some integer $m \geq 1$, $s$ is of the form
$$
s = e^{ic_1} \ldots e^{ic_m}, \ \text{ where } \ c_i = c_i^* \in A.
$$
This implies a self-adjoint analytic path $a(t) = a(t)^* \in C$ from $a_0$ to
$a_1$, for $t \in [0, 1]$, namely
$$
a(t) = e^{-ic_m t} \ldots e^{-ic_1 t} a_0 e^{ic_1 t} \ldots e^{ic_m t}.
$$
\endproclaim

\newpage

\proclaim{Corollary 13}
Let $\lambda_1, \dots, \lambda_n \in \Bbb R$.
Let $A$ be a unital complex $C^*$-algebra.
Let $a_0 \in S_p(A)$.
Then $a_0$ is in the centre of $A$ if and only if the connected component of
$a_0$ in $S_p(A)$ is $\{a_0\}$ (i.e., $a_0$ is isolated in $S_p(A)$).
\endproclaim

{\bf{Remark E.}}
For the $C^*$-algebra case, local connectedness of $S_p(A)$ via cubic
polynomial paths (analogously to Theorem 6 about $E_p(A)$) and connectedness
of connected components of $S_p(A)$ via polygonal paths (analogously to
Theorem~7 about $E_p(A)$) are false, already for $S(A)$.

Actually, this happens in the simplest case $A :=B(\Bbb C^2)$.
There the connected components of $S_p(A)$ are those with given rank, cf.\
[AMMZ03].
Thus one connected component is
$$
C:= \{T \in S(B(\Bbb C^2)) \mid \text{\rm rank } T = 1\}.
$$

For idempotent elements $T$, $\text{\rm rank }T = 1$ means that the
eigenvlaues are $0$, $1$ (with multiplicities $1$), or equivalently that the
sum and product of the eigenvalues are $1$, and $0$, respectively.
Equivalently, we have
$$
T_r T = 1 \ \text{ and } \ \text{\rm det } T = 0.
\tag{E.1}
$$

Self-adjointness of $T =:(a_{ij})_{i,j = 1}^2$ means $a_{11}, a_{22} \in \Bbb
R$, $a_{21} = \overline{a_{12}}$.
That is,
$$
\text{for some $a, b, c, d \in \Bbb R$ we have $a_{11} = a$, $a_{22} = b$,
$a_{12} = c + id$, $a_{21} = c - id$.}
$$

Rewriting (E.1) for $T$ self-adjoint, we have
$$
a + b = 1 \ \text{ and } \ ab - (c^2 + d^2) = 0.
\tag{E.2}
$$
Eliminating $b$, we obtain $a(1 - a) - (c^2 + d^2) = 0$, or, equivalently
$$
\frac{1}{4} = \left(a - \frac{1}{2}\right)^2 + c^2 + d^2.
\tag{E.3}
$$

Now suppose that for $t \in [0, 1]$ $T = T(t) =:(a_{ij}(t))_{i,j = 1}^2$ is a
polynomial of $t$.
Then, with the above notations, for $t \in [0, 1]$
$$
a = a(t), \ \ b = b(t), \ \ c = c(t), \ \ d = d(t)
$$
are polynomials of $t$ as well, and (E.3) holds for all $t \in [0, 1]$.
Now observe that if a polynomial of $t$ vanishes at more values of $t$ than
its degree, then it is identically~$0$.
This implies that equality (E.3) continues to hold for all $t \in \Bbb R$.
Then each of $a(t) - 1/2$, $a(t)$, $b(t)$, $(c(t)$, $d(t)$ is a polynomial
bounded on $\Bbb R$, i.e., is constant.

This shows that $C$ contains no non-constant polynomial path, while $|c| > 1$.
Hence $C$ is neither locally connected, nor connected via any polynomial paths.

\newpage

Moreover, we cannot prove Theorem 6 and the part of Theorem 7 concerning
similarities via finite products of exponentials for general Banach
${}^*$-algebras.
Namely, we used in their proofs S. Maeda [Ma], who already in the proof of his
Lemma~1, display, first equality uses the $C^*$-algebra property.
Probably the mentioned statements are false, but we have no counterexamples.

\vskip.1cm

The following Theorem 13 is a particular case of the next Theorem~14.
We separated these two statements because the more general Theorem~14 will be
proved by reducing its statement to its particular case dealt with in
Theorem~13.

\proclaim{Theorem 14}
Let $A = B(X)$, where $X$ is a complex Banach space.
Let $a_0, a_1 \in E_p(A)$, with $a_0$ fixed and $\|a_1 - a_0\|$ sufficiently
small.
Then there exists a polygonal path in $E_p(A)$, connecting $a_0$ and $a_1$,
consisting of $n$ segments.
The distance of this path from $a_0$ tends to $0$, if $\|a_1 - a_0\| \to 0$.

Hence the connected component of $a_0$ in $E_p(A)$ is locally pathwise
connected via paths of $n$ segments.
\endproclaim

\proclaim{Theorem 15}
Let $A$ be a unital complex Banach algebra.
Let $a_0, a_1 \in E_p(A)$, with $a_0$ fixed, and $\|a_1 - a_0\|$ small.
Then there exists a polygonal path in $E_p(A)$, connecting $a_0$ and $a_1$,
consisting of $n$ segments.
The distance of this path from $a_0$ tends to $0$, if $\|a_1 - a_0\| \to 0$.

Hence the connected component of $a_0$ in $E_p(A)$ is locally pathwise
connected via paths of $n$ segments.
\endproclaim

\proclaim{Theorem 16}
Let $A$ be a unital complex Banach algebra, let $C$ be a connected component
of $E_p(A)$, and let $a_0, a_1 \in C$.
Then there exists a polygonal path in $E_p(A)$, connecting $a_0$ and $a_1$.
\endproclaim

{\bf{Remark F.}}
For $C^*$-algebras, in general, polygonal connection between elements of a
connected component of $S_p(A)$ is impossible.
This holds already for $p(x) = x(x - 1)$, i.e., for $S(A)$, even in the
simplest case $A = B(H)$, for $H$ a Hilbert space.
Although this follows from Remark~E, we give here another argument, which
yields infinitely many examples for this special case.

\vskip.1cm

The connected components of $S(A) = S(B(H))$ are
$$
\align
&\bigl\{e \in S(A) \mid \text{\rm dim } N(E) = \alpha, \ \text{\rm dim } R(E)
= \beta,\\
&\quad \text{for some cardinalities } \alpha, \beta \geq 0, \text{ with
} \alpha + \beta = \text{\rm dim } H\bigr\}.
\endalign
$$
(Here dim is the dimension in Hilbert space sense.)
Let $\text{\rm dim }H \geq 2$, and let $0 < \beta = k < \text{\rm dim }H$ an
integer.
Let $C_k$ be the connected component of $S(A)$, consisting of self-adjoint
(i.e., orthogonal) projections of rank~$k$.
We claim that $C_k$ contains 

\newpage

no non-trivial segment (while it consists of more
than one element).

Let $a_0, a_1 \in C_k$, such that the segment $[a_0, a_1]$ lies in $C_k$.
For $a \in [a_0, a_1]$ $a$ is a compact operator, with singular values $\|a\|
= s_0(a) \geq s_1(a) \geq s_2(a) \geq \ldots$.
These $s_i(a)$'s are the eigenvalues of the non-negative square root
$\sqrt{a^* a} = a$, with multiplicities, in decreasing order.
(For $\text{\rm dim }H$ finite, $s_i(a) = 0$ for $i \geq \text{\rm dim }H$.)
Therefore $1 = s_0(a) = \ldots = s_{k - 1}(a) > 0 = s_k(a) = s_{k + 1}(a)
= \ldots$.
Thus all singular numbers are constant on the segment $[a_0, a_1]$.
This implies
$$
a_0 = a_1,
$$
by B. Aupetit, E. Makai, Jr., J. Zem\'anek [AMZ], Theorem, p.~517.

\proclaim{Corollary 17}
Let $A = B(H)$, where $H$ is a (complex) Hilbert space.
Then the (pathwise) connected components of $E_p(A)$ are of the form

\smallskip
{\narrower{\narrower{\narrower \noindent
$\{a \in E_p(A) \mid \text{ for each } 1 \leq i \leq n $
the Hilbert space dimension of the eigensubspace corresponding to the eigenvalue
 $\lambda_i$  is $\alpha_i\}$, where $\alpha_1, \dots, \alpha_n \geq 0$
  are any cardinalities whose sum is the Hilbert space dimension of~$H$.\par}}}

\smallskip
For the same $A$, considered as a $C^*$-algebra, and with
$\lambda_1, \dots, \lambda_n \in \Bbb R$, we have the following.
The (pathwise) connected components of $S_p(A)$ are of the same form as given
above, where still the respective eigensubspaces are orthogonal.

Therefore the connected components of $E_p(A)$ and $S_p(A)$ are just the
similarity classes of operators in $E_p(A)$ and $S_p(A)$ via conjugation with
an invertible, or unitary element.
\endproclaim

{\bf{Remark G.}}
For $A = B(H)$, if we consider all self-adjoint elements, with their
respective self-adjoint partitions of unity, then unitary similarity of these
self-adjoint partitions of unity for $a_0, a_1 \in A$ implies a path joining
$a_0$, $a_1$ in the set of elements unitarily similar to $a_0$, $a_1$.
This can be shown as in the second part of the proof of Corollary~16.
However it is not clear how to extend the definition of $S_p(A)$ from finite
real spectra to more general self-adjoint partitions of unity (and then to
investigate their connected components).

\vskip.1cm

{\bf{Problem.}}
It would be a natural conjecture that the distance of different connected
components of $E_p(A)$, for $A$ a unital complex Banach algebra (at least for
$A$ a unital complex $C^*$-algebra) and that the distance of different
connected components of  $S_p(A)$, for $A$ a unital complex $C^*$-algebra, and
$\lambda_1, \dots, \lambda_n \in \Bbb R$ is at least
$$
\min\bigl\{|\lambda_i - \lambda_j|\mid 1 \leq i, j \leq n, \ i \neq j\bigr\}.
$$

\newpage

This was the case for $p(\lambda) = \lambda(\lambda - 1)$, i.e., for
idempotents, or self-adjoint idempotents.
This implies the conjecture for $n = 2$, i.e., for $p(\lambda) = (\lambda
- \lambda_1)(\lambda - \lambda_2)$ in general.
We just have to observe that the set $E_{(\lambda - \lambda_2)(\lambda
- \lambda_2)}(A)$, or $S_{(\lambda - \lambda_n)(\lambda - \lambda_2)}(A)$ can
be obtained from $E(A)$, or $S(A)$, by the transformation
$x \mapsto \lambda_1 \cdot 1_A + (\lambda_2 - \lambda_1)x$, for $x \in A$.
For $S_{(\lambda - \lambda_1)(\lambda - \lambda_2)} (A)$, with
$\lambda_1, \lambda_2 \in \Bbb R$, we have by the same argument that the
distance of different connected components is at least $|\lambda_1
- \lambda_2|$.

Clearly this conjecture, if true, would be sharp, for any $A$.
Namely, $\lambda_1 \cdot 1_A, \dots, \lambda_n \cdot 1_A \in E_p(A)$ (or $\in
S_p(A)$), and since they are central, by Corollaries 8 and 12 their connected
components in $E_p(A)$ (in $S_p(A)$) are $\{\lambda_1 \cdot
1_A\}, \dots, \{\lambda_n \cdot 1_A\}$.
The minimal pairwise distance of these components is
$$
\min\bigl\{|\lambda_i - \lambda_j| \mid 1 \leq i,j \leq n, \ i \neq j\bigr\}.
$$

However, for $n \geq 3$ our proof does not give even that two different
connected components of $E_p(A)$ would have a positive distance.
There arise several questions.

1) Are these distances positive?

2) Are these distances bounded below by some positive function of
$\lambda_1, \dots, \lambda_n$?

3) Are these distances at least
$$
\min\bigl\{|\lambda_i - \lambda_j|\mid 1 \leq i,j \leq n, \ i \neq j\bigr\}?
$$
The same questions arise for the $C^*$-algebra case, with self-adjoint
idempotents, and $\lambda_1, \dots, \lambda_n \in \Bbb R$, but then questions
1) and 2) are answered positively in the following Theorem~16.

Even the commutative case would be of interest.

Observe that if we had a positive answer for 2) or 3), then this would imply
the stronger statement that even the spectral distance of different connected
components (i.e., the infimum of the spectral radii of the differences of
elements in the different components of $E_p(A)$ in question) would satisfy
the respective inequality.
(For the $C^*$-algebra case the norm equals the spectral radius, so there is
no such separate question.)
In fact, $E_p(A)$ and its connected components do not depend on the particular
norm of~$A$.
For $a \in C$, $b \in D$, with $C, D$ distinct components of $E_p(A)$ we can
renorm $A$ so that $\|a - b\| < \varrho(a - b) + \varepsilon$, for any
$\varepsilon > 0$ (and $\|1_A\|$ remains $1$), where $\varrho(\cdot)$ denotes
spectral radius.
Then any lower bound for $\| a - b\|$ implies the same lower bound for
$\varrho(a - b)$.

\vskip.1cm

If question 3) had a positive answer for Banach algebras, then it would have a
positive answer for  $C^*$-algebras as well.
Namely, different connected components $S_p(A)$ lie in different connected
components of $E_p(A)$, by [BFML], \S~1, Applications 2).
Observe that the considerations there only use the $C^*$-algebra generated by

\newpage

$\{P(t) \mid t \in [0, 1]\}$.
Also, they concern only the case of $E(A)$ and $S(A)$.
However, a path connecting $a_0, a_1 \in S_p(A)$ in $E_p(A)$ (cf.\ Theorem~7
and Theorem~16) yields paths connecting the idempotents $e_i(a_0)$ to the
idempotents $e_i(a_1)$, where $e_i(\cdot)$ are the functions in Proposition~1.
Moreover, by Corollary~5, $e_i(\cdot)$ from Proposition~1 is a continuous
function of its argument.

Since [BFML], \S~1, Applications 2) essentially uses the $C^*$-algebra
property, most probably for Banach ${}^*$-algebras with continuous involutions
different connected components of $S_p(A)$, or $S(A)$ may lie in the same
connected component of $E_p(A)$, or $E(A)$.
However, we do not have a concrete example.

\proclaim{Theorem 18}
Let $A$ be a unital complex $C^*$-algebra, and
$\lambda_1, \dots, \lambda_n \in \Bbb R$.
Then the distance of different connected components of $S_p(A)$ is at least
$$
\align
&\Bigl(\min_{1 \leq \alpha, \beta \leq n\atop \alpha \neq \beta}
|\lambda_\alpha - \lambda_\beta|\Bigr) \min_{1 \leq i \leq n} \prod_{j =
1 \atop j \neq i}^n |\lambda_i - \lambda_j| \Bigm/ \\
&\Bigm/ \biggl[\prod_{j = 1 \atop j \neq i}^n \Bigl(\max_{k \neq j} |\lambda_k
- \lambda_j | +
\min_{1 \leq \alpha, \beta \leq n \atop \alpha \neq \beta} |\lambda_\alpha
- \lambda_\beta|\Bigr) - \prod_{j = 1\atop j \neq i}^n \max_{k\neq j}
|\lambda_k - \lambda_j|\biggr].
\endalign
$$
\endproclaim

For $n = 2$ this gives back that the distance between different connected
components of $S(A)$ is at least $1$, and the distance between different
connected components of $S_p(A)$ is at least $|\lambda_1 - \lambda_2|$, that
is sharp,
cf.\ the following problem.
However, for $n \geq 3$ this estimate is probably far from the conjecturable
value
$$
\min_{i \neq j} |\lambda_i - \lambda_j|.
$$
E.g., for the case $\lambda_i = e^{(2\pi i)(i/n)}$ the expression in
Theorem~17 is
$$
e^{n ?? \left(\int\nolimits_0^\pi \log (2\sin t)dt/\pi + o(1)\right)} \bigm/
2^n,
$$
while
$$
\min_{i\neq j} |\lambda_i - \lambda_j| = |\lambda_2 - \lambda_1| =
2 \sin(\pi/n) = (2\pi / n) \cdot (1 + o(1)).
$$

Still observe that the estimate in Theorem~17 is invariant under simultaneous
translation of $\lambda_1, \dots, \lambda_n$, and their simultaneous
multiplication by a number of absolute value~$1$, and is homogeneous of first
degree under their simultaneous multiplication by a positive number.
These properties are shared by the (yet unknown) expression of the exact
infimum.

\newpage

Now we turn to another subject.
We begin with a notation.
For Banach spaces $X$, $Y$, we denote by $\Cal B(X, Y)$ the Banach space of
all bounded linear operators $X \to Y$.
In our papers B. Aupetit, E. Makai, Jr., M. Mbekhta, J. Zem\'anek [AMMZ03],
[AMMZ14] we investigated the following situation.
Let $A$ and $B$ be Banach algebras, and let $0 \in U \subset \Bbb C$ be an
open set.
We investigated analytic families of idempotents $q: U \to E(B) \subset B$,
and analytic familes $\pi: U \to \Cal B(A, B)$, whose values are surjective
Banach algebra homomorphisms $A \to B$.
We looked for conditions which assure that there exists an analytic family of
idempotents $p : U \to E(A) \subset A$, or $p : V \to E(A) \subset A$, where
$0 \in V \subset U$ is some open set of $\Bbb C$, such that
$$
q(\lambda) = \pi(\lambda) \ \text{ for all } \ \lambda \in U, \ \text{ or
} \ \lambda \in V, \text{ respectively}.
$$
(We called them global, or local liftings of analytic families of idempotents
along analytic families of surjective Banach algebra homomorphisms.)
We had analogous theorems for $S(A)$ and $S(B)$ as well.

It is a natural question, posed by L. W. Marcoux on the conference on linear
algebra in Ljubljana, Slovenia, 2014 (and also earlier by some participant of
a conference on operator theory at the Banach centre, Warsaw, some years ago),
whether these theorems have extensions for $E_p(A)$, or $S_p(A)$.

A {\it real analytic map} from an open subset $G$ of $\Bbb R$ $(\Bbb R^n)$ to
a Banach space is a map $f$ which for each $x_0 \in G$ has locally a power
series expansion
$$
f(x) = \sum_0^\infty a_i(x - x_0)^i
$$
(with the analogous formula for $\Bbb R^n$).
When we write in the following theorems spectrum of an element of $\text{\rm
Ker }\pi(\lambda) \subset A$ or $\text{\rm Ker }\pi(0) \subset A$, we mean
spectra in~$A$.
Two idempotents $e, f$ in some Banach algebra are {\it orthogonal} if $ef = fe
= 0$.
For examples where the hypotheses of the following theorems are satisfied cf.\
[AMMZ14], \S~3.

We remark that the global lifting Theorems 19, 20 are generalizations of\break 
[AMMZ14], Theorems 3, 4, if these are restricted to single idempotents, rather
than for sequences of idempotents.
However, the local lifting Theorems 21, 22 have stronger counterparts for
idempotents, cf.\ [AMMZ14], Theorems 1, 2, inasmuch there the spectral
hypotheses are weaker: the spectrum of each element of $\text{\rm Ker }\pi(0)$
does not disconnect $\Bbb C$, or is totally disconnected, respectively.

\proclaim{Theorem 19}
Let $U$ be an open subset of $\Bbb C$.
Let $A$ and $B$ be unital complex Banach algebras, and let $\pi: U \to \Cal
B(A, B)$ be an analytic map, whose  values are surjective unit-preserving
homomorphisms $A \to B$.
Suppose that the spectrum of  each element of $\text{\rm Ker }\pi(\lambda)$,
for each $\lambda \in U$, is $\{0\}$.
Let $b(\cdot) : U \to E_p(B)$ $(\subset B)$ be an analytic map.
Then there exists an analytic map $a(\cdot): U \to E_p(A)$ $(\subset A)$ such

\newpage

that
$$
\pi(\lambda) a(\lambda) = b(\lambda) \ \text{ for each } \ \lambda \in U.
$$
\endproclaim

\proclaim{Theorem 20}
Let $G$ be an open subset of~$\Bbb R$.
Let $A$ and $B$ be unital complex Banach algebras with continuous involutions,
and let $\pi : G \to \Cal B(A, B)$ be a real analytic map, whose values are
surjective unit-preserving ${}^*$-homomorphisms $A \to B$.
Suppose that the spectrum of each element of $\text{\rm Ker }\pi(\lambda)$,
for each $\lambda \in G$, is $\{0\}$.
Let $b (\cdot) : G \to S_p(B)$ $(\subset B)$ be a real analytic map.
Then there exists a real analytic map $a(\cdot) : G \to S_p(A)$ $(\subset A)$
such that
$$
\pi(\lambda) a(\lambda) = b(\lambda) \ \text{ for each } \ \lambda \in G.
$$
\endproclaim

The following Theorems 20 and 21 are localized versions of Theorems 18 and 19.
It is interesting that we have in Theorems 20 and 21 the spectral hypothesis
of Theorems 18 and 19 for $\lambda = 0$ only, and still we have local versions
of Theorems 18 and 19 for some open set ($V$ or $H$) containing~$0$.

\proclaim{Theorem 21}
Let $U$ be an open subset of $\Bbb C$, containing~$0$.
Let $A$ and $B$ be unital complex Banach algebras, and let $\pi : U \to \Cal
B(A, B)$ be an analytic map, whose values are unit-preserving homomorphisms
$A \to B$, such that $\pi(0)$ is surjective.
Suppose that the spectrum of each element of $\text{\rm Ker }\pi(0)$ is $\{0\}$.
Let $b(\cdot) : U \to E_p(B)$ $(\subset B)$ be an analytic map.
Then there exists an open set $V \subset \Bbb C$, such that $0 \in V \subset
U$, and an analytic map $a(\cdot) : V \to E_p(A)$ $(\subset A)$, such that
$$
\pi(\lambda) a(\lambda) = b(\lambda) \ \text{ for each } \ \lambda \in V.
$$
\endproclaim

\proclaim{Theorem 22}
Let $G$ be an open subset of $\Bbb R$, containing~$0$.
Let $A$ and $B$ be unital complex Banach algebras with continuous involutions,
and let $\pi: U \to \Cal B(A, B)$ be a real analytic map, whose values are
unit preserving ${}^*$-homomorphisms $A \to B$, such that $\pi(0)$ is
surjective.
Suppose that the spectrum of each element of $\text{\rm Ker }\pi(0)$ is $\{0\}$.
Let $b(\cdot) : G \to S_p(B)$ $(\subset B)$ be a real analytic map, then there
exist an open set $H \subset \Bbb R$, such that $0 \in H \subset G$, and a
real analytic map $a(\cdot) : H \to S_p(A)$ $(\subset A)$ such that
$$
\pi(\lambda) a(\lambda) = b(\lambda) \ \text{ for each } \ \lambda \in H.
$$
\endproclaim

{\bf{Remark H.}}
The following statements follow from the proofs in [AMMZ14], 

\newpage

cf.\ in
particular [AMMZ14], Remark~1.
In Theorem~18 we may replace $U$ by a Stein manifold.
In Theorem~19 we may suppose $G \subset \Bbb R^n$ open, provided each
connected component of $U$ has a neighbourhood base in $\Bbb C^n$ consisting
of domains of holomorphy, when $\Bbb R^n$ is embedded in $\Bbb C^n$ in the
canonical way.
In Theorem~20 we may suppose that $0 \in U \subset \Bbb C^n$ is open.
In Theorem~21 we may suppose that $0 \in G \subset \Bbb R^n$ is open.

\vskip.1cm

{\bf{Remark I.}}
We repeat a question posed in [\dots ] that has relations to the spectral
hypotheses in Theorems 19--22.
Let $H$ be a Hilbert space, and let $k(\lambda)$ be an analytic family of
compact operators in $\Cal B(H)$, for $\lambda \in \Bbb C$, with $|\lambda| <
1$.
Suppose that for some sequence $\{\lambda_n\} \subset \Bbb C$ with
$|\lambda_n| < 1$.
Converging to $0$, we have that the spectrum of $k(\lambda_n)$ is $\{0\}$.
Is then the spectrum of $k(\lambda)$ equal to $\{0\}$ {\it for each $\lambda$
with} $|\lambda| < 1$?
[\dots ] gave a positive answer if each $k(\lambda)$ has finite rank.

\vskip.1cm

{\bf{Remark J.}}
In [AMMZ03] we asked whether there exists an infinite dimensional Banach space
$X$, such that its Calkin algebra $\Cal C(X)$ (i.e., $\Cal B(X) / \Cal K(X)$,
where $\Cal B(X)$ and $\Cal K(X)$ are the Banach algebras of all bounded, or
compact linear operators on $X$, respectively) is commutative.
Then, in particular, by [Ze], cited in the Introduction under number 4, each
point of $E(\Cal C(X))$, being central, is isolated in $E(\Cal C(X))$.

It is known that there exists an infinite dimensional Banach space $X$, such
that
$$
\Cal B(X) = \{\lambda I + K \mid \lambda \in \Bbb C, \ K \in \Cal K(X)\},
$$
cf.\ [\dots].
Then $\Cal C(X) \cong \Bbb C$ is commutative, and thus $E(\Cal C(X)) = \{0,
1\}$ consists of two isolated points.

However, it remains open, whether $E(\Cal C(X))$ can consist of $n \in
[3, \infty)$ isolated points.

\heading
3. Proofs
\endheading

\demo{Proof of Proposition 1}
{\bf 1.} We begin with the proof of the implication that $a \in E_p(A)$
implies $a = \sum\limits_{i = 1}^n \lambda_i e_i$, with $\{e_1, \dots, e_n\}$
as in the proposition.

Let $p_1, \dots, p_n$ denote the Lagrange interpolation polynomials, i.e.,
$$
p_i(\lambda) = \prod_{j = 1\atop j \neq i}^n (\lambda
- \lambda_j) \Bigm/ \prod_{j = 1 \atop j \neq i}^n (\lambda_i - \lambda_j).
$$
Then

\newpage

$$
p :(\lambda_j) = \delta_{ij} \ \text{ for } \ 1 \leq i,j \leq n,
$$
and any polynomial $f$ over $F$, of degree at most $n - 1$ can be written, in
a unique way, as a linear combination of these polynomials, namely as
$$
f(\lambda) = \sum_{i = 1}^n f(\lambda_i) p_i(\lambda).
$$
In particular, we have
$$
1_F = \sum_{i = 1}^n p_i(\lambda), \ \text{ and } \ \lambda = \sum_{i =
1}^n \lambda_i p_i(\lambda).
\tag{1.1}
$$
(Recall that, by hypothesis, $n \geq 2$.)

Let us define
$$
e_i := p_i(a).
$$
Then $1_F = \sum\limits_{i = 1}^n p_i(\lambda)$ implies
$$
1_A = \sum_{i = 1}^n p_i(a) = \sum_{i = 1}^n e_i,
\tag{1.2}
$$
and $\lambda = \sum\limits_{i = 1}^n \lambda_i p_i(\lambda)$ implies
$$
a = \sum_{i = 1}^n \lambda_i p_i(a) = \sum_{i = 1}^n \lambda_i e_i.
\tag{1.2}
$$

Let $i \neq j$. Then $p_i(\lambda) p_j (\lambda)$ is a multiple of
$p(\lambda)$, and $p(a) = 0_A$, hence
$$
0_A = p_i(a) p_j(a) = e_i e_j .
\tag{1.4}
$$
Now we calculate $e_i^2$.
Let us divide $p_i(\lambda)$ by $\lambda - \lambda_i$ (with remainder),
obtaining
$$
p_i(\lambda) = (\lambda - \lambda_i) q_i(\lambda) + p_i(\lambda_i) = (\lambda
- \lambda_i) q_i(\lambda) + 1.
$$
Multiplying both sides with $p_i(\lambda)$, we obtain

\newpage

$$
p_i(\lambda)^2 = p_i(\lambda)(\lambda - \lambda_i) q_i(\lambda) + p_i(\lambda) =
\frac{p(\lambda)}{\prod\limits_{\scriptstyle i = 1\atop \scriptstyle j \neq
i}^n (\lambda - \lambda_j)} q_i(\lambda) + p_i(\lambda) .
$$
Substitute here $a$ for $\lambda$, and observe $p(a) = 0$.
Then we obtain
$$
e_i^2 = p_i(a)^2 = p_i(a) = e_i.
\tag{1.5}
$$

\smallskip
{\bf 2.} We turn to the proof of the converse implication.
So, let $a := \sum\limits_{i = 1}^n \lambda_i e_i$, with $\{e_1, \dots, e_n\}$
as in the proposition.
We calculate $p(a)$.
We have, using $\sum\limits_{i = 1}^n e_i = 1_A$, that
$$
\align
p(a) &= p \biggl(\sum_{i = 1}^n \lambda_i e_i\biggr) = \prod_{j =
1}^n \biggl(\sum_{i = 1}^n \lambda_i e_i - \lambda_j \biggr)
= \prod_{j = 1}^n \biggl(\sum_{i = 1}^n \lambda_i e_i - \lambda_j \cdot
1_A\biggr)\\
&= \prod_{j = 1}^n \biggl(\sum_{i = 1}^n \lambda_i e_i
- \lambda_j \cdot \sum_{i = 1}^n e_i \biggr) = \prod_{j = 1}^n \biggl(\sum_{i
= 1}^n (\lambda_i - \lambda_j) e_i \biggr).
\endalign
$$
Here the last expression is an $n$-fold product of sums of $n$ terms of the
form $(\lambda_i - \lambda_j)e_i$.
Rewrite this as a sum of $n^n$ terms, each term being an $n$-fold product, and
recall that the $e_i$'s commute.
Then each of these $n$-fold products, containing different $e_i$'s, vanishes.
There remains a sum of $n$ terms, each term being an $n$-fold product,
containing only a single~$e_i$.
Such an $n$-fold product equals
$$
\prod_{j = 1}^n (\lambda_i - \lambda_j) \cdot e_i^n = 0 \cdot e_i = 0_A.
$$
Then $p(a)$ is a sum of $n$ terms, each being equal to $0_A$, so
$$
p(a) = 0_A, \ \text{ i.e., } \ a \in E_p(A).
\tag{1.6}
$$

{\bf 3.} By {\bf 1} the function $a \mapsto e_i(a)$ is the polynomial function
$\lambda \mapsto p_i(\lambda)$.

\smallskip
{\bf 4.} By the hypothesis of the last statement of the proposition we have,
for each integer $N \geq 0$, that
$$
a^N = \sum_{i = 1}^n {\lambda_i}^N e_i = \sum_{j = 1}^m {\mu_j}^N
f_j, \ \text{ thus } \ 
\sum_{i = 1}^n {\lambda_i}^N e_i + \sum_{j = 1}^m {\mu_j}^N(-f_j) = 0.
$$

\newpage

We use these equations for $0 \leq N \leq n + m - 1$.
We have $n, m > 0$, and we may suppose $n \leq m$.

First suppose that $\lambda_i \neq \mu_j$ for each $1 \leq i \leq n$, $1 \leq
j \leq m$.
Then we have a system of linear equations on the vector space $A$, with
determinant the Vandermonde determinant of
$\lambda_1, \dots, \lambda_n, \mu_1, \dots, \mu_m$, that are all distinct.
Then the Vandermonde determinant is not $0$, and the unique solution of this
system is $e_1 = \ldots = e_n = f_1 = \ldots = f_m = 0$, a contradiction.
Hence, say, $\lambda_1 = \mu_1$.

Then our system of equations reduces to
$$
{\lambda_1}^N (e_1 - f_1) + \sum_{i = 2}^n {\lambda_i}^N e_i + \sum_{j = 2}^m
{\mu_j}^N (-f_j) = 0.
$$
If $\lambda_i \neq \mu_j$ for each $2 \leq i \leq n$, $2 \leq j \leq m$, then
we use these equations for $0 \leq N \leq n + m - 2$.
Again the determinant is the Vandermonde determinant of $\lambda_1$
$(=\mu_1)$, $\lambda_2, \dots, \lambda_n$, $\mu_2, \dots, \mu_m$.
Therefore
$$
e_1 - f_1 = e_2 = \ldots = e_n = f_2 = \ldots = f_m = 0.
$$
If $n = m = 1$, we have $\lambda_1 = \mu_1$, $e_2 = f_2$, and we are done.
Else we have a contradiction to $e_i \neq 0$, $f_j \neq 0$.
Therefore we may suppose, say, $\lambda_2 = \mu_2$.
We have
$$
{\lambda_2}^N(e_1 - f_1) + {\lambda_2}^N (e_2 - f_2) + \sum_{i = 3}^n
{\lambda_i}^N e_i + \sum_{j = 3}^m {\mu_j}^N f_j = 0.
$$
If $\lambda_i \neq \mu_j$ for each $3 \leq i \leq n$, $3 \leq j \leq m$, then
we use our equations only for $0 \leq N \leq n + m - 3$.
Like above, we obtain
$$
e_1 - f_1 = e_2 - f_2 = e_3 = \ldots = e_n = f_3 = \ldots = f_m = 0.
$$
If $n = m = 2$, we have $\lambda_1 = \mu_1$, $e_1 = e_2$, and also $\lambda_2
= \mu_2$, $e_2 = f_2$, and we are done.
Else we have a contradiction to $e_i \neq 0$, $f_j \neq 0$.
Therefore we may suppose, say, that $\lambda_3 = \mu_3$.

We continue analogously at each step, either the statement of the proposition
becomes proved, or we have a contradiction to $e_i \neq 0$, $f_j \neq 0$.

If we have not obtained the statement of the proposition earlier, we get to
$$
e_1 - f_1 = e_2 - f_2 = \ldots = e_n - f_n = f_{n + 1} = \ldots = f_m = 0.
$$
Then either $m = n$, and we have the statement of the proposition, or $m \geq
n + 1$, and we have a contradiction.
This proves the last statement of the proposition. \hfill $\blacksquare$
\enddemo

\newpage

\demo{Proof of Proposition 3}
In {\bf 1} of the proof of Proposition~1 we showed that $a \in E_p(A)$ implied
$a = \sum\limits_{i = 1}^n \lambda_i e_i$, with
$\{e_1, \dots, e_n\}$ as in Proposition~1.
Now we have $a \in S_p(A) \subset E_p(A)$, so these considerations apply.
There remains to prove ${e_i}^* = e_i$.
We have
$$
e_i = p_i(a) = \prod_{j = 1\atop j \neq i}^n (a - \lambda_j) \Bigm/ \prod_{j =
1 \atop j \neq i}^n (\lambda_i - \lambda_j).
$$
Therefore, $a = a^*$ and ${\lambda_i}^* = \lambda_i$ imply
$$
e_i = p_i(a) = p_i(a^*) = p_i(a)^* = {e_i}^*.
\tag{3.1}
$$

In {\bf 2} of the proof of Proposition~1 we showed that $a = \sum\limits_{i =
1}^n \lambda_i e_i$ with $\{e_1, \dots, e_n\}$ as in Proposition~1, implied
$a \in E_p(A)$.
Now $\lambda_i = {\lambda_i}^*$ and $e_i = {e_i}^*$ imply
$$
a^* = \Bigl(\sum_{i = 1}^n \lambda_i e_i\Bigr)^* = \sum_{i = 1}^n
{\lambda_i}^* {e_i}^* = \sum_{i = 1}^n \lambda_i e_i = a.
$$
Then $a \in E_p(A)$ and $a^* = a$ imply
$$
a \in S_p(A).
\tag{3.2}
$$

By Proposition~1 (cf.\ {\bf 3} of its proof), we know that $a \mapsto e_i(a)$
is a polynomial, with coefficients in~$F$.
But this polynomial is $\lambda \mapsto p_i(\lambda)$, whose coefficients are
rational functions of the self-adjoint elements $\lambda_i$, hence themselves
are self-adjoint.

\hfill $\blacksquare$
\enddemo

\demo{Proof of Theorem 6}
By Proposition~1 we have
$$
a_0 = \sum_{i = 1}^n \lambda_i e_{0i}, \ \ \ a_1 = \sum_{i = 1}^n \lambda_i
e_{1i},
$$
where $e_{0i}, e_{1i} \in E(A)$.
By Corollary~5, we have
$$
\|e_{1i} - e_{0i}\| < 1, \ \text{ for all } \ 1 \leq i \leq n,
$$
provided $\|a_1 - a_0\|$ is sufficiently small.
Therefore, by J. Zem\'anek [Ze], Lemma 3.1 and its proof, there exist $s_i \in
A$ involutions,
for all $1 \leq i \leq n$, such that

\newpage

$$
e_{1i} = s_i^{-1} e_{0i} s_i.
\tag{6.1}
$$
Moreover, for $\|a_1 - a_0\|$ sufficiently small we have $\|e_{1i} - e_{0i}\|$
sufficiently small, and then by [Ze], Lemma 3.1, Proof
$$
\|s_i - (2e_{0i} - 1)\| \ \text{ is sufficiently small as well.}
\tag{6.2}
$$
Here $2e_{0i} - 1$ is an involution, thus $\sigma(2e_{0i} - 1) \subset \{-1,
1\}$ so for $\|s_i - (2e_{0i} - 1) \|$ sufficiently small, by upper
semicontinuity of the spectrum,
$$
s_i = e^{c_i}, \ \text{ hence } \ e_{1i} := e^{-c_i} e_{0i} e^{c_i} ,
\tag{6.3}
$$
where $c_i \in A$.
Moreover, we may choose $c_i', c_i'' \in A$ so that
$$
e_{1i} = e^{-c_i''} e^{-c_i'} e_{0i} e^{c_i'} e^{c_i''} \text{ and } (c_i')^2
= (c_i'')^2 = 0,
\tag{6.4}
$$
cf.\ J. Esterle [Es], Theorem, Proof. Here even one of
the factors containing $c_i'$ can be deleted, without changing the value of
the right-hand side expression of the first equality in (6.4), leaving (6.4)
valid, cf.\ [Es], Theorem, Proof.
For $\|a_1 - a_0\|$ sufficiently small, we have for each $1 \leq i \leq n$
that $\|e_{1i} - e_{0i}\|$ is sufficiently small, and then for each $1 \leq
i \leq n$ we have that
$$
\text{both \ $\|c_i'\|$ \ and \ $\|c_i''\|$ \ are small,}
\tag{6.5}
$$
cf.\ [Es], Theorem, Proof.

Now we define
$$
s := \sum_{i = 1}^n e_{0i} s_i = \sum_{i = 1}^n s_i e_{1i} \ \text{ and } \ s'
:= \sum_{i = 1}^n s_i^{-1} e_{0i} = \sum_{i = 1}^n e_{1i} s_i^{-1}.
\tag{6.6}
$$
(The equalities (6.6) follow from (6.1).)
Then
$$
s's := \sum_{i = 1}^n s_i^{-1} e_{0i} s_i = \sum_{i = 1}^n e_{1i}  = 1, \ \ \
ss' := \sum_{i = 1}^n s_i e_{1i} s_i^{-1} = \sum_{i = 1}^n e_{0i} = 1,
\tag{6.7}
$$
hence
$$
s' = s^{-1}.
$$
Then we have

\newpage

$$
s^{-1} a_0 s = \sum_{i = 1}^n s_i^{-1} e_{0i} \cdot \sum_{j = 1}^n \lambda_j
e_{0j} \cdot \sum_{k = 1}^n e_{0k} s_k  =
\sum_{i = 1}^n s_i^{-1} e_{0i} \lambda_i s_i = \sum_{i = 1}^n \lambda_i e_{1i}
= a_1.
\tag{6.8}
$$
Rewriting this, by (6.6) and (6.3) we have
$$
a_1 = \sum_{i = 1}^n s_i^{-1} e_{0i} \cdot a_0 \cdot \sum_{i = 1}^n e_{0i} s_i
=
\sum_{i = 1}^n e^{-c_i} e_{0i} \cdot a_0 \cdot \sum_{i = 1}^n e_{0i} e^{c_i}.
\tag{6.9}
$$
Analogously, we have by (6.4)
$$
\aligned
a_1 &= \sum_{i = 1}^n \lambda_i e_{1i} = \sum_{i = 1}^n e^{-c_i''} e^{-c_i'}
e_{0i} \cdot a_0 \cdot \sum_{i = 1}^n e_{0i}
e^{c_i'} e^{c_i''}\\
&= \sum_{i = 1}^n (1 - c_i'') (1 - c_i') e_{0i} \cdot a_0 \cdot
\sum_{i = 1}^n e_{0i} \cdot (1 + c_i') (1 + c_i'').
\endaligned
\tag{6.10}
$$

This implies an analytic, or polynomial path
$$
a(t) := \sum_{i = 1}^n e^{-c_i t} e_{0i} \cdot a_0 \cdot \sum_{i = 1}^n e_{0i}
e^{c_it},
\tag{6.11}
$$
or
$$
\aligned
\tilde a(t) :&= \sum_{i = 1}^n e^{-c_i'' t}e^{-c_i t} e_{0i} \cdot
a_0 \cdot \sum_{i = 1}^n e_{0i} e^{c_i't} e^{c_i''t} \\
&=\sum_{i = 1}^n (1 - c_i''t) (1 - c_i't) e_{0i} \cdot a_0 \cdot \sum_{i =
1}^n e_{0i} (1 + c_i't)(1 + c_i'' t),
\endaligned
\tag{6.12}
$$
between $a(0) = a_0$ and $a(1) = a_1$, with $t \in [0, 1]$.
In (6.12), for all $1 \leq i \leq n$, one of the factors containing $c_i'$ can
be deleted, leaving the value of the last expression in (6.12) unchanged.
This means that we have a cubic polynomial path between $\tilde a(0) = a_0$
and $\tilde a(1) = a_1$.

Since by (6.5) $\|c_i'\|$ and $\|c_i''\|$ are small, for each $1 \leq i \leq
n$, the distance of this cubic polynomial path to $a_0$ for $t \in [0, 1]$ is
small, proving local connectedness of $E_p(A)$ via cubic polynomial paths.

For $\|a_1 - a_0\|$ sufficiently small, for all $1 \leq i \leq n$ \ $\|s_i -
(2e_{0i} - 1)\|$ is small (cf.\ (6.2)), and therefore

\newpage

$$
\align
\biggl\| \sum_{i = 1}^n e_{0i} \bigl(s_i - (2 e_{0i} - 1)\bigr)\biggr\|
&= \biggl\| \sum_{i = 1}^n e_{0i} s_i - \sum_{i = 1}^n e_{0i} \bigl(2 e_{0i} -
1\bigr) \biggr\| \\
&= \biggl\| \sum_{i = 1}^n e_{0i} s_i - \sum_{i = 1}^n e_{0i} \biggr\| = \|s -
1\|
\endalign
$$
is small as well.
As soon as $\|s - 1\| < 1$, we have
$$
s = e^c,
\tag{6.13}
$$
where $\| c\|$ is small as well, for $\|s - 1\|$ small.
Then, by (6.8) and (6.13)
$$
a_1 = s^{-1} a_0 s = e^{-c} a_0 e^c,
$$
that implies an analytic path
$$
\hat a(t) = e^{-ct} a_0 e^{ct}
\tag{6.14}
$$
between $\hat a(0) = a_0$ and $\hat a(1) = a_1$, for $t \in [0, 1]$.

For $\|a_1 - a_0\|$ sufficiently small, $\|c\|$ is sufficiently small, and
then the distance of this whole path from $a_0$ is small.
This proves local connectedness of $E_p(A)$, via similarities with exponential
functions.\hfill $\blacksquare$
\enddemo

{\bf{Remark K.}} 
Even though the inidividual $s_i$'s can be chosen as involutions, however most
probably
$$
s = \sum_{i = 1}^n e_{0i} s_i
$$
will not be an involution in general.
Thus in this respect $E_p(A)$ behaves differently from $E(A)$.
We do not know if $s$ could be chosen as an involution, in the proof of
Theorem~6.

The same remark concerns also the proof of Theorem 10, for the analogous
question for $C^*$-algebras.

\vskip.1cm

\demo{Proof of Theorem 7}
By Theorem 6, $E_p(A) \subset A$ is locally (pathwise) connected.
This implies that all connected components of $E_p(A)$ are relatively open
subsets of $E_p(A)$.

Now let us fix $a_0 \in C$.
We let
$$
C(a_0) := \bigl\{a\! \in\! E_p(A) \mid \exists m \geq 1 \text{ integer,
} \exists c_1, ..., c_m\! \in\! A, \
a = e^{-c_m}...  e^{-c_1} a e^{c_1} ... e^{c_m}\bigr\}.
\tag{7.1}
$$

\newpage

Since $a_0 = e^{-0} a_0 e^0$, we have $a_0 \in C$.
Then for $a \in C(a_0)$ there exists an analytic path
$$
a(t) = e^{-c_mt} \dots e^{-c_1 t} a_0 e^{c_1 t} \dots e^{c_m t},
\tag{7.2}
$$
hence $a \in C$.
That is,
$$
C(a_0) \subset C.
\tag{7.3}
$$

By Theorem 6, $C(a_0)$ is a relatively open subset of $E_p(A)$.
In fact, for $a \in C(a_0)$, any $a' \in E_p(A)$ sufficiently close to $a$ is
of the form
$$
a' = e^{-c} a e^c = e^{-c} e^{-c_m} \dots e^{-c_1} a_0 e^{c_1} \dots e^{c_m}
e^c,
$$
hence $a' \in C(a_0)$ as well (with integer $m + 1$ and elements $c_1, \dots,
c_m$, $c \in A$).

Now suppose that $C(a_0) \subsetneqq C$.
Since $C$ is connected, it cannot have a non-trivial open-and-closed subset.
Since $a_0 \in C(a_0) \subsetneqq C$, and $C(a_0)$ is open, therefore $C(a_0)$
is not closed.
That is, there is a point
$$
b \in \overline{C(a_0)} \cap \bigl(C \setminus C(a_0)\bigr)
\tag{7.4}
$$
($\overline{\phantom{CCC}}$ meaning closure in $E_p(A)$, that is closure in
$A$, by closedness of $E_p(A)$.)
By $b \in \overline{C(a_0)}$ any neighbourhood of $b$, relative to $E_p(A)$,
intersects $c(a_0)$.
By Theorem~6, some neighbourhood of $b$ consists of elements of the form
$e^{-c} b e^c$, where $c \in A$.
That is,
$$
e^{-c} b e^c = e^{-c_m} \ldots e^{-c_1} a_0 e^{c_1} \ldots e^{c_m},
\tag{7.5}
$$
which implies
$$
b = e^c e^{-c_m} \ldots e^{-c_1} a_0 e^{c_1} \ldots e^{c_m} e^{-c}.
\tag{7.6}
$$
Therefore $b \in C(a_0)$ (with integer $m + 1$, and elements $c_1, \dots, c_m,
-c \in A$).
However, by (7.4) $b \notin C(a_0)$, a contradiction.
This proves
$$
C(a_0) = C.
\tag{7.7}
$$

Then for $a_1 \in C = C(a_0)$ (7.1) implies an analytic path $a(t) \in
E_p(A)$, thus $a(t) \in C$, from $a_0$ to $a_1$.
Namely,
$$
a(t) = e^{-c_m t} \ldots e^{-c_1 t} a_0 e^{c_1 t} \ldots e^{c_m t},
\tag{7.8}
$$
for $t \in [0,1]$.

As we have seen in the proof of Theorem~6, (6.4), (6.12), we may additionally
suppose $c_1^2 = \ldots = c_m^2 = 0$, which implies a polynomial path

\newpage

$$
\tilde a(t) = e^{-c_mt}\ldots e^{-c_1 t} a_0 e^{c_1 t} \ldots e^{c_m t} = (1 -
c_m t) \ldots (1 - c_1 t) a_0(1 + c_1 t) \ldots (1 + c_m t),
\tag{7.9}
$$
from $a_0$ to $a_1$.
For $t \in [0, 1]$, and here one of the factors containing $c_1$ can be
deleted, without changing the value of the right-hand side of (7.9), by
J. Esterle [Es], Theorem, Proof (cf.\ the proof of Theorem~6).\hfill
$\blacksquare$
\enddemo

\demo{Proof of Corollary 8}
Let $a_0 \in E_p(A)$, and let $C$ be the connected component of $a_0$ in
$E_p(A)$.

Suppose that $a_0$ is  in the centre of $A$, then for any $a_1 \in C$, by
Theorem~7, for some invertible $s$ we have
$$
a_1 = s^{-1} a_0 s = a_0 s^{-1} s = a_0.
$$
Hence $C = \{a_0\}$.

Now suppose that $a_0$ is not in the centre of~$A$.
Let $x \in A$, $a_0 x \neq x a_0$.
Then, for $t \to 0$,
$$
C \ni e^{-tx} a_0 e^{tx} = (1 - tx) a_0 (1 + tx) + O(t^2) = a_0 + t (-xa_0 +
a_0 x) + O(t^2) \neq a_0.
$$
Then $C$ is not the singleton $\{a_0\}$.\hfill $\blacksquare$
\enddemo

\demo{Proof of Theorem 9}
Let $a \in E_p(A)$ be an arbitrary element of $E_p(A)$. That is,
$$
a = \sum_{i = 1}^n \lambda_i e_i,
$$
where $\{e_1, e_2, e_3, \dots, e_n\}$ is an arbitrary partition of unity to
mutually orthogonal idempotents in~$A$.
We want to perturb $e_1$ and $e_2$ to $e_1'$ and $e_2'$.
Leaving $e_3, \dots, e_n$ unchanged, so that we obtain a new partition of
unity to mutually orthogonal idempotents in~$A$.

Let $x \in A$ be arbitrary, and let
$$
e_1' := e_1 + e_1 x e_2, \ \ e_2' := e_2 - e_1 x e_2, \ \text{ and } \ e_i' =
e_i \text{ for } i \geq 3.
$$
Then
$$
e_1' + e_2' + e_3' + \ldots + e_n' = 1, \ \ (e_1')^2 = e_1', \ \ (e_2')^2 =
e_2', \ \ e_1'e_2' = e_2'e_1' = 0,
$$
and for $i,j \geq 3$, $i \neq j$
$$
e_1' e_i' = e_i'e_1' = e_2' e_i' = e_i'e_2' = 0 \ \text{ and } \ e_i'e_j' = 0.
$$

\newpage

Thus $\{e_1', e_2', e_3', \dots, e_n'\}$ is a new partition of unity to
mutually orthogonal idempotents in~$A$.

We consider the element
$$
a' := \lambda_1 e_1' + \lambda_2 e_2' + \lambda_3 e_3' + \ldots + \lambda_n
e_n' \in E_p(A).
$$
Multiplying it from left by $e_1', e_2'$ we obtain
$$
e_1' a' = \lambda_1 e_1' = \lambda_1(e_1 + e_1 x e_2) \ \text{ and } \ e_2' a'
= \lambda_2 e_2' = \lambda_2(e_2 - e_1 x e_2).
$$
Also we have either $\lambda_1 \neq 0$ or $\lambda_2 \neq 0$.
In the first case we use the equation for $e_1'a'$, in the second case we use
the equation for $e_2'a'$.
Thus unless the continuous linear map $A \ni x \mapsto \lambda_1 \cdot e_1 x
e_2$ (or $\lambda_2 \cdot e_1 x e_2$) is identially $0$, its image contains a
(complex) line.
In the second one of these cases also the image of the continuous linear map
$$
\align
A \ni x \mapsto a' &= \lambda_1 e_1' + \lambda_2 e_2' + \lambda_3 e_3'
+ \ldots + \lambda_n e_n'\\
&= \lambda_1(e_1 + e_1 x e_2) + \lambda_2(e_2 - e_1 x e_2) + \lambda_3 e_3
+ \ldots + \lambda_n e_n
\endalign
$$
contains a (complex) line.
For $x = 0$ we have $a'= \lambda_1 e_1 + \lambda_2 e_2 + \lambda_3 e_3
+ \ldots + \lambda_n e_n$, so the above line contans the arbitrarily chosen
element $a$ of $E_p(A)$.
Now $a'\in E_p(A)$, and $a'$ is in the same connected component of $E_p(A)$ as
$a$.

We can replace $e_1, e_2$ by any other $e_i, e_j$ $(i \neq j)$ in the above
considerations.
Therefore either

1) for each $x \in A$ and each $i \neq j$ we have $e_i x e_j = 0$, or

2) for some $i \neq j$ the image of the map $A \ni x \mapsto \sum\limits_{i =
1}^n \lambda_i e_i'$ contains a complex line.

In Case 2) the theorem is proved.

In Case 1)
$$
\align
e_i x &= e_i x \cdot 1 = \sum_{i \neq j}^n e_i x e_j + e_i x e_i = e_i x
e_i \ \text{ and }\\
x e_i &= 1 \cdot x e_i = \sum_{i \neq j}^n e_j x e_i + e_i x e_i = e_i x e_i.
\endalign
$$
Therefore
$$
e_i x = e_i x e_i = x e_i.
$$

\newpage

Here $x \in A$ is arbitrary, hence $e_i \in Z(A)$, and
$$
a = \sum_{i = 1}^n \lambda_i a_i \in Z(A).
$$
Then by Corollary 8 the connected component of $E_p(A)$ containing $a$ is the
singleton~$\{a\}$.\hfill $\blacksquare$
\enddemo

\demo{Proof of Corollary 10}
By Corollary 8 and Theorem 9, the central points of $E_p(A)$ are isolated in
$E_p(A)$, and the non-central points of $E_p(A)$ lie on complex lines
contained in $E_p(A)$. \hfill $\blacksquare$
\enddemo

\demo{Proof of Theorem 11}
The proof is analogous to that of Theorem~6.

We have that $s_i$, that was an involution in the proof of Theorem~6, is now
even a self-adjoint involution, by S. Maeda [Ma], Lemma 2.
Moreover, $s$, defined in (6.6) in the proof of Theorem~6, now satisfies, also
using Proposition~3, (6.7) from the proof of Theorem~6, and $s_i^2 = 1$, that
$$
s^* = \sum_{i = 1}^n s_i^* e_{0i}^* = \sum_{i = 1}^n s_i e_{0i} = \sum_{i =
1}^n s_i^{-1} e_{0i} = s' = s^{-1},
\tag{11.1}
$$
so $s$ is unitary as well.
Moreover, from the proof of Theorem~6, for $\| a_1 - a_0\|$ sufficiently
small, $\|s - 1\|$ is small as well.
As soon as $\|s - 1\| < 1$, we have
$$
s = e^{ic},
\tag{11.2}
$$
where $\|c\|$ is small for $\|s - 1\|$ small. 
In fact, letting $\log$ the branch of the logarithm function that takes the
value $0$ at $1$, we define
$$
ic := \log s.
\tag{11.3}
$$
By (10.1)
$$
c = c^*.
\tag{11.4}
$$

Then
$$
a(t) = a(t)^* = e^{-ict} a_0 e^{ict}
\tag{11.5}
$$
is, for $t \in [0, 1]$, a self-adjoint analytic path between $a_0$ and $a_1$.

Like in Theorem~6, for $\|a_1 - a_0\|$ sufficiently small, the distance of
this whole path from $a_0$ is small, showing local connectedness of $S_p(A)$,
by similarities via unitary valued exponential functions.\hfill $\blacksquare$
\enddemo

\demo{Proof of Theorem 12}
Like in the deduction of the global Theorem~7 from the local 

\newpage

Theorem~6, we
have from the local Theorem~9 that
$$
a_1 = e^{-ic_m} \ldots e^{-ic_1} a_0 e^{ic_1} \ldots e^{ic_m},
$$
where $c_i = c_i^* \in A$ for $2 \leq i \leq m$.

The statement about the self-adjoint analytic path follows immediately.\hfill
$\blacksquare$
\enddemo

\demo{Proof of Corollary 13}
Let $a_0 \in S_p(A)$ and let $C$ be the connected component of $a_0$ in
$S_p(A)$.

Suppose that $a_0$ is in the centre of $A$.
Then, like in the proof of Corollary~8 we obtain $C = \{a_0\}$.
(By the way, this actually follows from Corollary~8.)

Now suppose that $a_0$ is not in the centre of~$A$.
Let $x \in A$, $a_0 x \neq x a_0$.
Let $x = y + i_2$, $y = y^*$, $z = z^*$.
Then $a_0$ does not commute either with $y$, or with $z$.
That is, we can assume $x = x^*$.
Then, like in the proof of Corollary~8 for $t \to 0$ we have
$$
C \ni e^{-itx} a_0 e^{itx} = a_0 + it(-x a_0 + a_0 x) + O(t^2) \neq a_0.
$$
Then $C$ is not the singleton $\{a_0\}$.\hfill $\blacksquare$
\enddemo

\demo{Proof of Theorem 14}
{\bf 1.} We have by Proposition 1 and Corollary 5
$$
a_0 = \sum_{i = 1}^n \lambda_i e_{0i}, \ \  a_1 = \sum_{i = 1}^n \lambda_i
e_{1i}, \text{ with $\|e_{1i} - e_{0i}\|$ small for each } 1 \leq i \leq n.
\tag{14.1}
$$
By Corollary~2, we have
$$
X = \bigoplus_{i = 1}^n X_{0i} = \bigoplus_{i = 1}^n X_{1i},
\tag{14.2}
$$
with coordinate projections
$$
e_{0i} : X \to X_{0i}, \ \ \ e_{i1} = X \to X_{1i}.
$$
We will consider $n + 1$ decompositions of $X$ as direct sums, namely
$$
\aligned
X &= \bigoplus_{i = 1}^n X_{0i}, \ \ X = X_{11} \oplus \biggl(\bigoplus_{i =
2}^n X_{0i}\biggr),\\
X &= X_{11} \oplus X_{12} \oplus \biggl( \bigoplus_{i = 3}^n
X_{0i}\biggr), \dots, \\
X &= \biggl( \bigoplus_{i = 1}^{n - 1} X_{1i}\biggr) \oplus X_{0n}, \ \ X
= \bigoplus_{i = 1}^n X_{1i}.
\endaligned
\tag{14.3}
$$

\newpage

We will show that these are in fact direct decompositions of $X$, and that all
their respective coordinate projections are close to $e_{01}, \ldots, e_{0n}$.
Moreover, we will show that for any two neighbouring decompositions in (13.3),
say, the $j^{\text{\rm th}}$  and $(j + 1)^{\text{\rm st}}$  one, where
$j \in \{1, \ldots, n\}$, a simultaneous linear interpolation between their
first, \dots, $n^{\text{\rm th}}$  coordinate projections, for a parameter
$t \in [0, 1]$, gives a linear path of $n$ pairwise orthogonal idempotents,
summing to $I \in B(X)$, for each $t \in [0, 1]$.

Hence, multiplying these $n$ orthogonal idempotents by
$\lambda_1,\ldots, \lambda_n$, respectively, and summing these, we obtain a
linear path in $E_p(A)$, joining its values for $t = 0$ and $t = 1$.
These values are obtained from the first, \dots, $n^{\text{\rm th}}$
projections of the $j^{\text{\rm th}}$  and $(j + 1)^{\text{\rm st}}$
decompositions (for some $1 \leq j \leq n$), by multiplying them by
$\lambda_1, \dots, \lambda_n$, respectively, and summing them.

We will consider the case $j = 1$, i.e.\ the construction of the second direct
decomposition of $X$ from the first one in (14.3).
The further cases are analogous.

{\bf 2.} Let $e_0, e_1 \in E(A) = E(B(X))$, where $e_0$ is fixed, and $\|e_1 -
e_0\|$ is small.
Then the Kovarik element
$$
g = g(e_0, e_1)
\tag{14.4}
$$
is defined in the following way, cf.\ Z. V. Kovarik [Ko], p.~343, (7),
Definition of $S_1, C_1$ and p.~345, (b), and p.~347, Definition of $F$, and
J. Esterle [Es], p.~253, Theorem, Proof, first paragraph.
Suppose $\|e_1 - e_0\| < 1$.
Then $e_0 + e_1 - 1$ is invertible ([Ko], p.~345, (b) and [Es], above cited),
and we let
$$
g = g(e_0, e_1) := e_1(e_0 + e_1 - 1)^{-2} e_0 = e_1\bigl[1 - (e_1 -
e_0)^2\bigr]^{-1} e_0,
\tag{14.5}
$$
then
$$
\aligned
g &= g^2 \in E(A),\\
e_1 g &= g, \ \text{ thus } \ R(g) \subset R(e_1), \\
g e_1 &= e_1, \ \text{ thus } \ R(e_1) \subset R(g), \ \text{ hence } \ R(g) =
R(e_1), \\
e_0 g &= e_0, \ \text{ thus } \ N(g) \subset N(e_0),\\
g e_0 &= g, \ \text{ thus } \ N(e_0) \subset N(g), \ \text{ hence } \ N(g) =
N(e_0)
\endaligned
\tag{14.6}
$$
(here $R(\cdot)$ is range, $N(\cdot)$ is kernel), cf.\ [Ko], p.~347, [Es],
p.~253).
Clearly, a projection $g$ is uniquely determined by $R(g) = R(e_1)$ and $N(g)
= N(e_0)$
for each $t \in [0, 1]$ (even for each $t \in \Bbb C$) we have
$$
(1 - t) e_0 + tg, \ (1 - t) e_1 + tg \in E(A)
\tag{14.7}
$$
cf.\ [Ko], p.~347.
(By the way, all these calculations are rather straightforward.)

{\bf 3.} As told in {\bf 1}, we want to change the direct sum decomposition $X
= \bigoplus\limits_{i = 1}^n X_{0i}$ to the direct sum decomposition
$X = X_{10} \oplus \Bigl(\bigoplus\limits_{i = 2}^n X_{0i}\Bigr)$, in such a
way that the 

\newpage

simultaneous linear interpolation (for $t \in [0, 1]$) of the
first, \dots, $n^{\text{\rm th}}$  projections gives for all $t \in [0, 1]$ a
system of orthogonal projections summing to~$1$.

The projections corresponding to the direct sum decomposition $X
= \bigoplus\limits_{i = 1}^n X_i$ are $e_{01}, \dots, e_{0n}$.
The projection to $X_{11}$, corresponding to the direct sum decomposition
$$
X = X_{11} \oplus \biggl(\bigoplus_{i = 2}^n X_{0i}\biggr),
\tag{14.8}
$$
has range $X_{11} = R(e_{11})$, and has kernel $\bigoplus\limits_{i = 2}^n
X_{0i} = N(e_{01})$, hence supposing $\|e_{11} - e_{01}\| < 1$, it exists and
equals
$$
g(e_{01}, e_{11}) = e_{11} (e_{01} + e_{11} - 1)^{-2} e_{01}.
\tag{14.9}
$$
Since later we will use only $g(e_{01}, e_{11})$, we most often will write for
it just~$g$.
Thus
$$
R(g) = X_{11} \ \text{ and } \ N(g) = \bigoplus_{i = 2}^n X_{0i}.
\tag{14.10}
$$
Observe that for $\|a_1 - a_0\|$ small we have $\|e_{1i} - e_{0i}\|$ small for
each $1 \leq i \leq n$, and $\|g(e_{01}, e_{11}) - e_{01}\|$ small.

For the complementary projection $1 - g$ we have
$$
R(1 - g) = \bigoplus\limits_{i = 1}^n X_{0i} \ \text{ and } \ N(1 - g) = X_{11}.
\tag{14.11}
$$
Now we consider
$$
e_{0i} (1 - g) \ \text{ for } \ 2 \leq i \leq n.
\tag{14.12}
$$
This is obtained by first applying the projection $1 - g$, having range
$\bigoplus\limits_{i = 2}^n X_{0i}$, and then applying the restriction to
$\bigoplus\limits_{i = 2}^n X_{0i}$ of the $i^{\text{rm th}}$ projection $X
= \bigoplus\limits_{i = 1}^n X_{0i} \to X_{0i}$.
This shows that $e_{0i}(1 - g)$ is a projection of $X$ to $X_{0i}$, so
$$
R\bigl(e_{0i} (1 - g)\bigr) = X_{0i}.
\tag{14.13}
$$
Similarly one sees that

\newpage

$$
X_{11} = N(1 - g) \subset N\bigl(e_{0i}(1 - g)\bigr) \text{ and for } 2 \leq
j \leq n, \ j \neq i \text{ we have } X_{0j} \subset N\bigl(e_{0i} (1 -
g)\bigr).
$$
Hence
$$
N\bigl( e_{0i} (1 - g)\bigr) = X_{11} \oplus\biggl( \bigoplus_{j = 2 \atop
j \neq i}^n X_{0j}\biggr).
\tag{14.14}
$$

Now we turn the opposite way.
We define $g(e_{01}, e_{11})$ by (14.10), which is possible as soon as
$\|e_{11} - e_{01}\| < 1$.
Then we have
$$
R(g) = R(e_{11}) \ \text{ and } \ N(g) = (e_{01}).
\tag{14.15}
$$
We will show that
$$
g \ \text{ and } \ e_{0i}(1 - g) \ \text{ for } \ 2 \leq i \leq n
\tag{14.16}
$$
form an orthogonal system of idempotents summing to~$1$.
Their sum is
$$
g + \sum_{i = 2}^n e_{0i} (1 - g) = g + (1 - e_{01}) (1 - g) = 1 + e_{01} g -
e_{01} = 1
\tag{14.17}
$$
by (13.6).

{\bf 4.} It remained to show that the elements in (13.16) form an orthogonal
system of idempotents.
We will show more: a simultaneous linear interpolation between $e_{01}$ and
$g$ and between $e_{0i}$ and $e_{0i}(1 - g)$, for $t \in [0, 1]$, gives an
orthogonal system of idempotents, summing to~$1$.

That is
$$
(1 - t) e_{01} + tg \text{ and } (1 - t) e_{0i} + te_{0i} (1 - g) \text{ for }
2 \leq i \leq n
\tag{14.18}
$$
are orthogonal systems of idempotents summing to~$1$.

Since $e_{01} + \sum\limits_{i = 2}^n e_{0i} = 1$ and $g + \sum_{i = 2}^n
e_{0i}(1 - g) = 1$ (cf.\ (14.17)), the sum of all elements in (14.16) is $1$
for all $t \in [0, 1]$.

The inclusion
$$
(1 - t) e_{01} + tg(e_{01}, e_{11}) \in E(A)
\tag{14.19}
$$
is contained in Z. V. Kovarik [Ko], p.~347.

For $2 \leq i \leq n$ we have
$$
(1 - t) e_{0i} + te_{0i} (1 - g) = e_{0i} - t e_{0i} g.
\tag{14.20}
$$
Its square is

\newpage

$$
e_{0i} - t e_{0i} g - t e_{0i} g e_{0i} + t^2 e_{0i} g e_{0i} g.
$$
Here the third and fourth summands contain
$$
g(e_{01}, e_{11}) e_{0i} = e_{11} (e_{01} + e_{11} - 1)^{-2} e_{01} e_{02} = 0.
$$
Therefore the element (14.20) belongs to~$E(A)$.

There remains to show that the elements in (14.18) are pairwise orthogonal.

Let $2 \leq i \leq n$.
Then, using the shorter form in (14.19),
$$
[(1 - t)e_{01} + tg ] \cdot [e_{0i} - te_{0i}g]
= (1 - t) e_{01} e_{0i} - (1 - t) t e_{01} e_{0i} g + tg e_{02} - t^2 g e_{02}
g.
\tag{14.21}
$$
Here the first and second summand contain $e_{01} e_{0i} = 0$, and the third
and fourth summands contain $g e_{02} = 0$.
Hence (14.21) is~$0$.

On the other hand,
$$
\aligned
&[e_{0i} - te_{0i} g]\cdot [(1 - t)e_{01} + tg] \\
&= (1 - t) e_{0i} e_{01} + te_{0i} g - t(1 - t) e_{0i} g e_{01} - t^2 e_{0i} g
g \\
&= 0 + t e_{0i} g - t(1 - t) e_{0i} g - t^2 e_{0i} g = e_{0i} g \cdot (t - t(1
- t) - t^2) = 0.
\endaligned
\tag{14.22}
$$

Now let $2 \leq i, j \leq n$, $i\neq j$.
Then
$$
\aligned
&[e_{0i} - te_{0i}g] \cdot [e_{0j} - te_{0j} g]\\
&= e_{0i} e_{0j} - te_{0i} e_{0j} g - t e_{0i}g e_{0j} + t^2 e_{0i} g e_{0j} g.
\endaligned
\tag{14.23}
$$
Here the first and second summands contain $e_{0i} e_{0j} = 0$, and the third
and fourth summands contain
$$
g(e_{01}, e_{11}) e_{0j} = e_{11} (e_{01} + e_{11} - 1)^2 e_{01} e_{0j} = 0.
$$

{\bf 5.}
Then multiplying the simultaneous interpolating elements in (12.15)??
(14.15)??, which form for each $t \in [0,1]$ a partition of unity, by the
respective $\lambda_i$'s, and summing them, we obtain a segment in $E_p(A)$,
beginning from $a_0$.
We have to show that its other endpoint is
$\sum\limits_{i = 1}^n \lambda_i f_i$, where the $f_i$'s are the projections
associated to the direct sum decomposition $X =
X_{11} \oplus \Bigl( \bigoplus\limits_{i = 2}^n X_{0i}\Bigr)$.
For $i = 1$ we have seen this in {\bf 3}, (14.10), while for $2 \leq i \leq n$
we have seen this in {\bf 3}, (14.13) and (14.14).
This means that our procedure changed the direct sum representation $X
= \bigoplus\limits_{i = 1}^n X_{0i}$ to 

\newpage

$X =
X_{11} \oplus \Bigl(\bigoplus\limits_{i = 2}^n X_{0i}\Bigr)$, as we wanted to
show.
Moreover, the respective projections changed just a bit.
Then, for $\|a_1 - a_0\|$ sufficiently small, we can change in the second
analogous step $X_{02}$ to $X_{12}$, \dots, in the $n^{\text{\rm th}}$
analogous step $X_{0n}$ to $X_{1n}$, the last direct sum decomposition will be
thus $X = \bigoplus\limits_{i = 1}^n X_{1i}$.
Even in this last decomposition the respective projections remain close to the
original projections $e_{0i}$, showing local connectedness of $E_p(A)$ by
polygons consisting of $n$ segments.
Therefore, we obtained in $E_p(A)$ a polygon of $n$ segments, joining $a_0
= \sum\limits_{i = 1}^n \lambda_i e_0$ to $a_1 = \sum\limits_{i =
1}^n \lambda_i f_i$, where $\{e_1, \dots, e_n\}$, $\{f_1, \dots, f_n\} \subset
E(A)$ are partitions of unity.

{\bf 6.} As mentioned in {\bf 5}, by the first exchange all the respective
projections changed just a bit, so they remained close to $e_{01}, \dots,
e_{0n}$.
The same holds for the second, \dots, $n^{\text{\rm th}}$ exchange.
Then linearly interpolating between the successive respective first, \dots,
$n^{\text{\rm th}}$ projections, we obtain polygonal paths close to the
original projections $e_{01}, \dots, e_{0n}$.\hfill $\blacksquare$
\enddemo

\demo{Proof of Theorem 15}
{\bf 1.}
We use the same formulas as in the proof of Theorem~13.
That is, we define from $a_0 = \sum\limits_{i = 1}^n \lambda_i e_{0i}$ and
$a_1 = \sum\limits_{i = 1}^n \lambda_i e_{1i}$, where $\|a_1 - a_0\|$ is
small, hence by Corollary~5 also $\|e_{1i} - e_{0i}\|$ for each $1 \leq i \leq
n$ is small, the elements
$$
\gathered
g = g(e_{01}, e_{11}),\ \ 1 - g, \ \ e_{0i}(1 - g), \ \ (1 - t) e_{01} +
tg, \ \ (1 - t) e_{0i} + t e_{0i}(1 - g)\\
\text{ for each } 2 \leq i \leq n, \text{ and each } t \in [0,1].
 \endgathered
\tag{15.1}
$$

Recall that in the proof of Theorem~12 we exchanged the direct summands
$X_{0i}$ with the direct summands $X_{1i}$, one by one, in the first step for
$i = 1, \dots$, for the $n^{\text{\rm th}}$ step for $i = n$.
However, we did not directly use these subspaces, but only the respective
projections, and the whole proof was done by algebraic manipulations with
these projections.
At the end of the proof it was established that after the last step of these
$n$ exchanges $X_{01}, \dots, X_{0n}$ was exchanged to $X_{11}, \dots,
X_{1n}$.
Equivalently, the first, \dots, $n^{\text{\rm th}}$ coordinate projection for
the direct sum decomposition $X = \bigoplus\limits_{i = 1}^n X_{0i}$ were
exchanged, after the last of the $n$ steps, by the first, \dots, $n^{\text{\rm
th}}$ coordinate projection for the direct sum decomposition $X
= \bigoplus\limits_{i = 1}^n X_{1i}$.

We have to show that performing these exchanges, the $n$-tuple
$e_{01}, \ldots, e_{0n}$ will change just to $e_{11}, \ldots, e_{1n}$, also
for the case of a Banach algebra, and not only for the case of the algebra of
bounded operators on a Banach space.

\newpage

We recapitulate the formulas describing these exchanges in a form suitable for
us.
We consider an $(n + 1) \times n$ matrix, with entries in $E(A)$, namely
$$
\left(
\matrix
f_{01} & \ldots & f_{0n}\\
f_{11} & \ldots & f_{1n}\\
\vdots & & \vdots\\
f_{n1} & \ldots & f_{nn}
\endmatrix
\right).
\tag{15.2}
$$
We set $f_{01} = e_{01}, \dots, f_{0n} = e_{0n}$.
For $1 \leq i \leq n$ we set
$$
\aligned
f_{ii} &= g(f_{i - 1, i}, e_{1i}) \ \text{ for } \ 1 \leq i \leq n, \ \text{
and }\\
f_{ij} &= f_{i - 1, j} (1 - g(f_{i - 1, i}, e_{1i})) \ \text{ for } \ 1 \leq
i,j \leq n, \ \ i \neq j.
\endaligned
\tag{15.3}
$$

Observe that for $e_{01}, \ldots, e_{0n}$ fixed and $\|e_{1j} - e_{0j}\|$
small we have for each $1 \leq i,j \leq n$ that $\|f_{ij} - e_{0j}\|$ is
small.
This follows by induction for $i$.
For $j = i$ this follows by continuity of the function $(e, f) \mapsto g(e,f)$
(where $e, f \in E(A)$ are close).
For $j \neq i$, by the same continuity $f_{ij}$ is close to $e_{0j}(1 -
g(e_{0i}, e_{0i})) = e_{0j} (1 - e_{0i}) = e_{0j}$.

Using the recursive formulas (14.3), we see that each $f_{ij}$ $(0 \leq i \leq
n$, $1 \leq j \leq n)$ is a ``rational'' function of $e_{01}, \ldots, e_{0n},
e_{11}, \ldots, e_{1n}$.
More exactly, each $f_{ij}$ can be expressed by $e_{01}, \dots, e_{0n},
e_{11}, \dots, e_{1n}$ by using the following operations:
linear combinations, multiplication, and applying the function $x \mapsto (1 -
x)^{-1} = \sum\limits_{m = 0}^\infty x^m$, where $\|x\|$ is small
($\|x\| < 1$ is necessary).
Therefore in particular $f_{nj}$, for each $1 \leq i \leq n$, belongs to the
closed subalgebra of $A$, generated by $e_{01},\dots, e_{0n}, e_{11}, \dots,
e_{1n}$.
We have to show $f_{nj} = e_{1j}$ for each $2 \leq j \leq n$.
(Theoretically there would be a possibility to express $f_{21},\dots, f_{2n}$,
then further $f_{31}, \dots, f_{3n}$, etc.\
by these ``rational'' functions, superposed to each other, but these formulas
soon become very complicated, and seem not to be treatable.)

{\bf 2.} Let us consider the regular representation
$$
\varphi: \ A \to B(A),
\tag{15.4}
$$
that maps any $a \in A$ to the (bounded) operator $x \mapsto ax$ (i.e., left
multiplication by $a$).
This is a Banach algebra homomorphism, and even is an isometric embedding
$A \to B(A)$.
In particular, $\varphi$ is injective.

Then $\varphi$ preserves linear combinations, multiplication and the function
$x \mapsto (1 - x)^{-1}$ (for $\|x\|$ small).
Therefore considering $\varphi e_{01}, \dots, \varphi e_{0n}, \varphi
e_{11}, \dots, \varphi e_{1n} \in E(B(A)) \subset B(A)$, and also
$\varphi f_{ij} \in E(B(A)) \subset B(A)$, for each $0 \leq i \leq n_n$ and
$1 \leq j \leq n$, we have the following.
Each $\varphi f_{ij}$, $0 \leq i \leq n$, $1 \leq j \leq n$, thus in
particular each $\varphi f_{nj}$, $1 \leq j \leq n$, is expressed by the same
formula via $\varphi e_{01}, \dots, \varphi e_{0n}$, 

\newpage

$\varphi
e_{11}, \dots, \varphi e_{1n}$, as $f_{ij}$, in particular $f_{nj}$, is
expressed via $e_{01}, \dots, e_{0n}$, $e_{11}, \dots, e_{1n}$.
Observe that $\varphi e_{01}, \dots, \varphi e_{0n}, \varphi
e_{11}, \dots, \varphi e_{1n}$ and
$\varphi f_{n1},\dots, \varphi f_{nn}$ can be considered as to play the role
of $e_{01}, \dots, e_{0n}, e_{11}, \dots, e_{1n}$ and
$f_{n1}, \dots, f_{nn}$ for the case of the operator algebra $B(A)$, rather
than for~$A$.

However, by Theorem~14 we already know the equalities $f_{nj} = e_{1j}$ for
$1 \leq j \leq n$, for the case of $A = B(X)$, where $X$ is a Banach space.
In particular, this holds for $B(A)$, i.e., we have
$$
\varphi f_{nj} = \varphi e_{1j} \ \text{ for } \ 1 \leq j \leq n.
\tag{15.5}
$$
Since $\varphi$ is injective, this implies
$$
f_{nj} = e_{1j} \ \text{ for } \ 1 \leq j \leq n,
\tag{15.6}
$$
that was to be shown.

{\bf 3.}
It remained to prove that the linear interpolations between $f_{0i}, f_{1i}$,
and between $f_{1i}, f_{2_i}$, \dots, and between $f_{n - 1,i}, f_{ni}$
together give a polygonal path which?? is close to the original projection
$f_{0i} = e_{0i}$, for each $1 \leq i \leq n$.
However, this was proved in Theorem 13 for $B(X)$ for any Banach space X, thus
in particular for $B(A)$.
Now recall that the regular representation $\varphi : \ A \to B(A)$ is an
isometry into.
From Theorem~13 we know that the image by $\varphi$ of this polygonal path in
$B(A)$ is close to $\varphi e_{0i}$, for each $1 \leq i \leq n$.
By the isometric property of $\varphi$ this polygonal path in $B$ is close to
$E_{0i}$, for each $1 \leq i \leq n$.\hfill $\blacksquare$
\enddemo

\demo{Proof of Theorem 16}
The proof follows from Theorem 14 in an analogous way as the proof of
Theorem~7 followed from Theorem~6 (and as the proof of Theorem~11 followed
from Theorem~10).\hfill $\blacksquare$
\enddemo

\demo{Proof of Corollary 17}
Let $C$ be a connected component of $E_p(A)$, or of $S_p(A)$, and let $a_0 \in
C$.
Then for any $a_1 \in C$, $a_1$ is similar to $a_0$ via some invertible, or
unitary element of $A = B(H)$, cf.\ Theorems 7 and 11.
This similarity implies the equality of the Hilbert space dimensions of the
eigensubspaces of $a_0$ and $a_1$ corresponding to the eigenvalue $\lambda_i$,
for each $1 \leq i \leq n$.

Conversely, let for $a_0, a_1 \in E_p(A)$, or $a_0, a_1 \in S_p(A)$ (in the
second case with $\lambda_1, \dots, \lambda_n \in \Bbb R$) the Hilbert space
dimensions of the eigensubspaces corresponding to each eigenvalue $\lambda_i$,
$1 \leq i \leq n$, be equal.
Then $H$ is the direct sum of these eigensubspaces, both for $a_0$ and $a_1$.
In the $C^*$-algebra case these eigensubspaces are even orthogonal, then we have
$$
a_1 = s^{-1} a_0 s,
$$
where $s \in B(H)$ is invertible, and in the $C^*$-algebra case $s$ is unitary.
Now recall that $Ge(H)$, as well as the unitary group $U(H)$ of $H$, is
(pathwise) connected.

\newpage

Therefore we have a path in $Ge(H)$, or $U(H)$, from $Id$ to $s$.
Its image by the continuous map $t \mapsto t^{-1} a_0t$ is a path from $a_0$
to $a_1$, in $E_p(A)$, or $S_p(A)$.
For the $C^*$-algebra case (in the second case with
$\lambda_1, \dots, \lambda_n \in \Bbb R$).\hfill $\blacksquare$
\enddemo

\demo{Proof of Theorem 18}
We proceed analogously to Remark C, 2), only will make use of the
$C^*$-algebra property.
Also we recall the proof of Theorem 6, and use the notations there.

So we have, also using Propositions 1 and 3 that
$$
a_0 = \sum_{i = 1}^n \lambda_i e_{0i} \sum_{i = 1}^n \lambda_i e_i(a_0), \ \
a_1 = \sum_{i = 1}^n \lambda_i e_{1i} = \sum_{i = 1}^n \lambda_i e_i(a_1).
\tag{18.1}
$$
As soon as
$$
\|e_i(a_1) - e_i(a_0)\| < 1, \ \text{ for each } \ 1 \leq i \leq n,
\tag{18.2}
$$
we have that $e_i(a_0), e_i(a_1)$ belong to the same connected component of
$S(A)$,
and their linear combinations
$$
a_0 = \sum_{i = 1}^n \lambda_i e_i (a_0), \ \  \ a_1 = \sum_{i =
1}^n \lambda_i e_i (a_1)
$$
belong to the same connected component of $S_p(A)$.
(Observe that then $(e_1(a_0),\dots, e_n(a_0))$ and $(e_1(a_1), \dots,
e_1(a_n))$ belong to a connected component of $S(A)^n$, as the product of
connected sets is connected.
Then its image by the linear, thus continuous function $(e_1, \dots,
e_n) \mapsto \sum\limits_{i = 1}^n \lambda_i e_i$ is also connected.)

We will write
$$
a_1 = a_0 + \min_{1 \leq \alpha, \beta \leq n\atop \alpha \neq \beta}
|\lambda_\alpha - \lambda_\beta| \cdot \varepsilon,
\tag{18.3}
$$
where by the problem before Theorem~17 we may suppose $\|\varepsilon\| < 1$.
We look for upper estimates of $\|\varepsilon\|$ which assure (18.2). ??
Recall that
$$
e_i(\lambda) = \prod_{j = 1 \atop j \neq i}^n (\lambda
- \lambda_j) \Bigm/ \prod_{j = 1\atop j \neq i}^n (\lambda_i - \lambda_j).
\tag{18.4}
$$
Hence for $a \in A$, $a = \sum\limits_{k = 1}^n \lambda_k e_k$ we have by
$\sum\limits_{k = 1}^n e_k = 1$ that the numerator of $e_i(a)$ is
$$
\prod_{j = 1\atop j \neq i}^n (\lambda - \lambda_j) =
\prod_{j = 1\atop j \neq i}^n \biggl(\sum_{k = 1}^n \lambda_k e_k
- \lambda_j\biggr) =
\prod_{j = 1\atop j \neq i}^n \biggl(\sum_{k = 1}^n (\lambda_k - \lambda_j)
e_k\biggr).
\tag{18.5}
$$

\newpage

Hence the numerator of $\|e_i(a_0 + \varepsilon) - e_i(a_0)\|$ is
$$
\aligned
&\biggl\| \prod_{j = 1\atop j \neq i}^n \Bigl(a - \lambda_j
+ \min_{1 \leq \alpha, \beta \leq n \atop \alpha \neq \beta} |\lambda_\alpha
- \lambda_\beta | \cdot \varepsilon \Bigr) -
\prod_{j = 1\atop j \neq i}^n (\alpha - \lambda_j)\biggr\| \\
&= \prod_{j = 1 \atop j \neq i}^n \biggl( \sum_{k = 1}^n (\lambda_k
- \lambda_j ) e_k + \min_{1 \leq \alpha, \beta \leq n\atop \alpha \neq \beta}
|\lambda_\alpha - \lambda_\beta| \cdot \varepsilon\biggr) - \prod_{j = 1\atop
j \neq i}^n \sum_{k = 1}^n (\lambda_k - \lambda_j) e_k\biggr\|.
\endaligned
\tag{18.6}
$$
(18.6) is a non-commutative polynomial of $\varepsilon$, but $\sum\limits_{k =
1}^n (\lambda_k - \lambda_j) e_k$ are self-adjoint commuting elements.
Clearly in (18.6) ?? the coefficient of $\varepsilon^\circ$ is~$0$.
All other terms contain $\varepsilon$ at least once.
Performing the multiplications in
$$
\prod_{j = 1\atop j\neq i}^n \biggl(\sum_{k = 1}^n (\lambda_k - \lambda_j) e_k
+ \min_{1 \leq \alpha , \beta \leq n\atop \alpha \neq \beta} |\lambda_\alpha
- \lambda_\beta|\varepsilon \biggr),
\tag{18.7}
$$
we obtain $2^{n - 1}$ additive terms, each of degree $n - 1$.
Each of these additive terms is a product or $n - 1$ factors, each factor
being of the form
$$
\sum_{k = 1}^n (\lambda_k - \lambda_j) e_k,  \ \text{ or } \
\min_{1 \leq \alpha, \beta \leq n \atop \alpha \neq \beta} |\lambda_\alpha
- \lambda_\beta|\varepsilon.
\tag{18.8}
$$
The norm of a single such product is at most the product of the norms of the
factors.
This last product contains some power of
$\bigl\|\min\limits_{1 \leq \alpha, \beta\leq n \atop \alpha \neq \beta}
|\lambda_\alpha - \lambda_\beta| \varepsilon\bigr\|$, say,
$$
\bigl\|\min_{1 \leq \alpha, \beta \leq n\atop \alpha \neq \beta}
|\lambda_\alpha - \lambda_\beta| \varepsilon \bigr\|^k,
$$
which we estimate from above by
$$
\Bigl(\min_{1 \leq \alpha, \beta \leq n \atop \alpha \neq \beta}
|\lambda_\alpha - \lambda_\beta|\Bigr)^k - \|\varepsilon\|
\tag{18.9}
$$
(observe that now $k \geq 1$, and $\|\varepsilon \| < 1$, so
$\|\varepsilon \|^k \leq \| \varepsilon\|$).
The upper estimate of the last mentioned single product still contains $n - 1
- k$ factors of the form
$$
\biggl\| \sum_{k = 1}^n (\lambda_k - \lambda_j) e_k\biggr\| = \max_{1 \leq
k \leq n} |\lambda_k - \lambda_j|,
\tag{18.10}
$$

\newpage

(i.e., here $j$ can take $n - 1 - k$ values),
the equality holding by the $C^*$-algebra property and $e^*_k = e_k$ (cf.\
Proposition~3).
(We note yet that if at least two factors of the form $\sum_{k = 1}^n
(\lambda_k - \lambda_j) e_k$ follow each other in a summand after performing
the multiplications in (18.7), then further simplifications become possible:
since $e_k^2 = e_k$ and $e_{k_1} e_{k_2} = 0$ for $k_1 \neq k_2$, so such a
product will be a linear combination of the $e_k$'s only.
However, such a better result would give great complications in the estimate
of Theorem~18.)

Then (18.6) can be estimated above
$$
\aligned
&\biggl\|\prod_{j = 1\atop j\neq i}^n \biggl(\sum_{k = 1}^n (\lambda_k
- \lambda_j) e_k + \min_{1 \leq \alpha, \beta \leq n\atop \alpha \neq \beta}
|\lambda_\alpha - \lambda_\beta| \biggr) - \prod_{j = 1 \atop j \neq
i}^n \sum_{k = 1}^n (\lambda_k - \lambda_j)
e_k \biggr\| \cdot \|\varepsilon \| \\
&\leq \biggl[\prod_{j = 1 \atop j\neq i}^n \biggl(\Bigl\| \sum_{k = 1}^n
(\lambda_k - \lambda_j)e_k \Bigr\| + \min_{1 \leq \alpha, \beta \leq
n\atop \alpha \neq \beta} |\lambda_\alpha - \lambda_\beta|\biggr) - \prod_{j =
1\atop j \neq i}^n \Bigl\| \sum_{k = 1}^n (\lambda_k - \lambda_j)
e_k \Bigr\| \biggr] \cdot \|\varepsilon \|.
\endaligned
\tag{18.11}
$$
(Observe that here writing the first product as the sum of $2^{n - 1}$
summands, one summand cancels with the second product, which corresponds to
the fact that in (18.7) the coefficient of $\varepsilon^\circ$ is~$0$.)
Then (18.11) with (18.10) imply
$$
\biggl[\prod_{j = 1\atop j\neq i}^n \biggl(\max_{1 \leq \alpha, \beta \leq
n\atop \alpha \neq \beta} |\lambda_k - \lambda_j|
+ \min_{1 \leq \alpha, \beta \leq n\atop \alpha \neq \beta} |\lambda_\alpha
- \lambda_\beta|\biggr) - \prod_{j = 1\atop j \neq i}^n \max_{1 \leq k \leq n}
|\lambda_k - \lambda_j| \biggr] \cdot \|\varepsilon \| \biggm/ \prod_{j =
1 \atop j \neq i}^n |\lambda_i - \lambda_j|.
\tag{18.12}
$$

Now recall (18.12) and the considerations following it.
Then, as soon as
$$
\|e_i(a_1) - e_i(a_0)\| = \|e_i(a_0 + \varepsilon) - e_i (a_0)\| < 1 \ \text{
for each } \ 1 \leq i \leq n,
\tag{18.13}
$$
$a_0$ and $a_1$ belong to the same connected component of $S_p(A)$.
That is, for
$$
\max_{1 \leq i \leq n} \|e_i(a_1) - e_i (a_0)\| < 1
$$
$a_0$ and $a_1$ belong to the same connected component of $S_p(A)$.
Said otherwise, if $a_0$ and $a_1$ belong to different connected components of
$S_p(A)$, then
$$
\max_{1 \leq i \leq n} \|e_i (a_i) - e_i (a_0)\| \geq 1.
$$

\newpage

Then (18.13) gives for $\min\limits_{1 \leq \alpha, \beta \leq
n\atop \alpha \neq \beta} |\lambda_\alpha - \lambda_\beta|$ \
$\|\varepsilon \| = \|a_1 - a_0\|$ the lower estimate given in the
theorem.\hfill $\blacksquare$
\enddemo

We give the proof of Theorem~19, and then will indicate the necessary changes
to prove Theorems 20, 21 and 22.

\demo{Proof of Theorem 19}
Let $\lambda \in U$, thus $b(\lambda) \in E_p(B)$.
By Proposition~1, ``only if'' part
$$
\gather
b(\lambda) = \sum_{i = 1}^n \lambda_i f_i(\lambda) = \sum_{i = 1}^n \lambda_i
e_i(b(\lambda)),\\
\text{where }
\{f_1(\lambda), \dots, f_n(\lambda)\} \subset E(B) \text{ is a partition of
unity},
\endgather
$$
with $e_i(\cdot)$ the polynomials from Proposition~1.
That is, $f_1(\lambda), \dots, f_n(\lambda)$ are mutually orthogonal, and
$f_1(\lambda) + \ldots + f_n(\lambda) = 1$.
Observe that $U \ni \lambda \mapsto e_i(b(\lambda)) \in E(B) \subset B$, as a
composition of two analytic functions, is analytic.
By [AMMZ14], Theorem~3, using the hypotheses of this theorem, this analytic
family of mutually orthogonal idempotents admits a lifting $\tilde
f_i(\cdot): \ U \to E(A) \subset A$, which is also an analytic family of
mutually orthogonal idempotents, thus satisfies
$$
\pi(\lambda) \tilde f_i(\lambda) = f_i(\lambda) \ \text{ for each
} \ \lambda \in U.
$$
Of course we cannot guarantee $\tilde f_1(\lambda) + \ldots + \tilde
f_n(\lambda) = 1$.
Therefore we replace $\tilde f_n(\cdot)$ by $1 - \tilde f_1(\cdot) - \ldots
- \tilde f_{n - 1}(\cdot)$, which is also an analytic family of idempotents,
and which forms with $\tilde f_1(\cdot), \dots, \tilde f_{n - 1}(\cdot)$
a mutually orthogonal system of idempotents, summing to~$1$, i.e., we have a
partition of unity in $E(A)$.

Once more we apply Proposition~1, now the ``if'' part.
Thus
$$
\lambda_1 \tilde f_1(\cdot) + \ldots + \lambda_{n - 1} \tilde f_{n - 1}(\cdot)
+ \lambda_n (1 - \tilde f_1(\cdot) - \ldots - \tilde f_{n - 1} (cdot)) \in
E_p(A),
$$
and thus we have in the last display an analytic family of elements of $E_p(A)$.
Still we have to prove that this analytic family lifts $b(\cdot)$.
In fact,
$$
\align
&\pi(\lambda) \bigl(\lambda_1 \tilde f_1(\lambda) + \ldots + \lambda_{n -
1} \tilde f_{n - 1}(\lambda)\bigr)
+ \lambda_n(1 - \tilde f_1(\lambda) - \ldots - \tilde f_{n = 1}(\lambda)\bigr)\\
&= \lambda_1 f_1(\lambda) + \ldots + \lambda_{n - 1} f_{n - 1}(\lambda)
+ \lambda_n\bigl(1 - f_1(\lambda) - \ldots - f_{n - 1}(\lambda)\bigr) \\
&= \lambda_1 f_1(\lambda) + \ldots + \lambda_{n - 1} f_{n - 1}(\lambda)
+ \lambda_n f_n(\lambda) = b(\lambda).
\endalign
$$
\hfill $\blacksquare$
\enddemo

\demo{Proof of Theorem 20}
We replace in the proof of Theorem~19 Proposition~1 by Proposition~3, and
[AMMZ14] Theorem~3 by [AMMZ14], Theorem~4.\hfill $\blacksquare$
\enddemo

\newpage

\demo{Proof of Theorem 21}
We use Proposition~1 and [AMMZ14], Theorem~5.
Thus we obtain liftings $\tilde f_1(\lambda), \dots, \tilde f_{n - 1}
(\lambda)$, $1 - \tilde f_1(\lambda) - \ldots - \tilde f_{n - 1}(\lambda)$,
not for each $\lambda \in U$, but only for $\lambda \in
V_1, \dots, \lambda \in  V_{n - 1}$, $\lambda \in V_1 \cap \ldots \cap V_{n -
1}$, respectively, for some open subsets $V_1, \dots, V_{n - 1}$ of $U$, each
containing~$0$.
Then Theorem~21 holds for $V := V_1 \cap \ldots \cap V_{n - 1}$.\hfill
$\blacksquare$
\enddemo

\demo{Proof of Theorem 22}
We use Proposition~3, and [AMMZ14], Theorem~6, and finish like in the proof of
Theorem~21.\hfill $\blacksquare$
\enddemo

\heading{References}\endheading

\parskip=2.5pt

\eightpoint
\frenchspacing
\leftskip=12mm
\item{\hbox to14mm{[Au]\hfill}} B. Aupetit, Projections in real Banach algebras,
{\it Bull. London Math. Soc.} {\bf 13} (1981), 412--414.

\item{\hbox to14mm{[ALZ]\hfill}} B. Aupetit, T. J. Laffey, J. Zem\'anek,
Spectral classification of projections,
{\it Linear Alg. Appl.} {\bf 41} (1981), 131--135.

\item{\hbox to14mm{[AMZ]\hfill}}  B. Aupetit, E. Makai, Jr., J. Zem\'anek,
Strict convexity of the singular value sequence,
{\it Acta Sci. Math. (Szeged)} {\bf 62} (1996), 517--521.

\item{\hbox to14mm{[AMMZ03]\hfill}} B. Aupetit, E. Makai, Jr., M. Mbekhta,
J. Zem\'anek,
The connected components of the idempotents in the Calkin algebra, and their
liftings, operator theory and Banach algebras,
{\it Conf. Proc., Rabat (Morocco), April 1999} (Eds. M. Chidami, R. Curto,
M. Mbekhta, F.-H. Vasilescu, J. Zem\'anek), Theta, Bucharest, 2003, 23--30.

\item{\hbox to14mm{[AMMZ14]\hfill}} B. Aupetit, E. Makai, Jr., M. Mbekhta,
J. Zem\'anek,
Local and global liftings of analytic families of idempotents in Banach
algebras,
{\it Acta Sci. Math. (Szeged)} {\bf 80} (2014), 149--174.

\item{\hbox to14mm{[AZ]\hfill}} B. Aupetit, J. Zem\'anek,
On zeros of analytic multivalued functions,
{\it Acta Sci. Math. (Szeged)} {\bf 46} (1983), 311--316.

\item{\hbox to14mm{[BFML]\hfill}} Z. Boulmaarouf, M. Fernandez Miranda,
J.-Ph. Labrousse,
An algorithmic approach to orthogonal projections and Moore--Penrose inverses,
{\it Numer. Funct. Anal. Optim.} {\bf 18} (1997) (1-2), 55--63, MR97m:65105.

\item{\hbox to14mm{[Ca]\hfill}} J. W. Calkin,
Two-sided ideals and congruences in the ring of bounded operators in Hilbert
space,
{\it Ann. of Math.} {\bf 42} (1941), 839--873.

\item{\hbox to14mm{[Dao]\hfill}} A. Daoui,
Sur le degr\'e minimum des connections polynomiales entre les projections dans
une alg\`ebre de Banach,
{\it Rend. Circ. Mat. Palermo} {\bf 49} (2000), 353--362.

\item{\hbox to14mm{[Dav]\hfill}} C. Davis,
Separation of two linear subspaces,
{\it Acta Sci. Math. (Szeged)} {\bf 19} (1958), 172--187.

\item{\hbox to14mm{[dlH]\hfill}} P. de la Harpe,
Initiation \`a l'alg\`ebre de Calkin, in:
{\it Alg\`ebres d'op\'erateurs, Les Plans-sur-Bex, 1978}, ed. P. de la Harpe,
Lecture Notes in Math. {\bf 725}, Springer, Berlin, 1979, 180--219.

\item{\hbox to14mm{[En]\hfill}} {\it Encyclopaedia of Math.}, Vol. 2.
An updated and annotated translation of the {\it Soviet Mathematical
Encyclopaedia}, Editor-in-chief I. M. Vinogradov, Kluwer, Dordrecht etc.,
1988.

\item{\hbox to14mm{[Es]\hfill}} J. Esterle,
Polynomial connections between projections in Banach algebras,
{\it Bull. London Math. Soc.} {\bf 15} (1983), 253--254.

\newpage

\item{\hbox to14mm{[EG]\hfill}} J. Esterle, J. Giol,
Polynomial and polygonal connections between idempotents in finite-dimensional
real algebras,
{\it Bull. London Math. Soc.} {\bf 36} (2004), 378--382.

\item{\hbox to14mm{[HLZ]\hfill}} V. K. Harchenko, T. J. Laffey, J. Zem\'anek,
A characterization of central idempotents,
{\it Bull Acad. Polon. Sci., S\'er. Sci. Math.} {\bf 29} (1981), 43--46.

\item{\hbox to14mm{[Ho91]\hfill}} J. P. Holmes,
Idempotents in differentiable semigroups,
{\it J. Math. Anal. Appl.} {\bf 162} (1991), 255--267.

\item{\hbox to14mm{[Ho92]\hfill}} J. P. Holmes,
The structure of the set of idempotents in a Banach algebra,
{\it Illinois J. Math.} {\bf 36} (1992), 102--115.

\item{\hbox to14mm{[La]\hfill}} J.-Ph. Labrousse,
The general local form of an analytic mapping into the set of idempotent
elements of a Banach algebra,
{\it Proc. Amer. Math. Soc.} {\bf 123} (1995), 3467--3471.

\item{\hbox to14mm{[IKa]\hfill}} I. Kaplansky,
{\it Fields and rings}, The University of Chicago Press, Chicago, 1972.

\item{\hbox to14mm{[TKa]\hfill}} T. Kato,
{\it Perturbation theory for linear operators},
Springer, Berlin, 1966.

\item{\hbox to14mm{[YKa75]\hfill}} Y. Kato,
An elementary proof of Sz. Nagy's theorem,
{\it Math. Japon.} {\bf 20} (1975), 257--258.

\item{\hbox to14mm{[YKa76]\hfill}} Y. Kato,
Some theorems on projections of von Neumann algebras,
{\it Math. Japon.} {\bf 21} (1976), 367--370.

\item{\hbox to14mm{[Ko]\hfill}} Z. V. Kovarik,
Similarity and interpolation between projectors,
{\it Acta Sci. Math. (Szeged)} {\bf 39} (1977), 341--351.

\item{\hbox to14mm{[Ku]\hfill}} K. Kuratowski,
{\it Topology II}, Academic Press, New York--London; PWN Warsaw, 1968.

\item{\hbox to14mm{[Ma]\hfill}} S. Maeda,
On arcs in the space of projections of a $C^*$-algebra,
{\it Math. Japon.} {\bf 21} (1976), 371--374.

\item{\hbox to14mm{[MZ]\hfill}} E. Makai, Jr., J. Zem\'anek,
Nice connecting paths in connected components of sets of algebraic elements in
a Banach algebra,
{\it Czechoslovak J. Math.}, accepted.

\item{\hbox to14mm{[MZ89]\hfill}} E. Makai, Jr., J. Zem\'anek,
On polynomial connections between projections,
{\it Lin. Alg. Appl.} {\bf 126} (1989), 91--94.

\item{\hbox to14mm{[Ra]\hfill}} I. Raeburn,
The relation between a commutative Banach algebra and its maximal ideal space,
{\it J. Funct. Anal.} {\bf 25} (1977), 366--390.

\item{\hbox to14mm{[Ri]\hfill}}  C. E. Rickart,
{\it General theory of Banach algebras},
van Nostrand, New York, 1960.

\item{\hbox to14mm{[RSzN]\hfill}} F. Riesz, B. Sz\H{o}kefalvi Nagy,
{\it Le\c cons d'analyse fonctionnelle},
Akad\'emiai Kiad\'o, Budapest, 1955; III-\`eme \'ed., Paris, Gauthier-Villard.

\item{\hbox to14mm{[SzN42]\hfill}} B. Sz\H{o}kefalvi Nagy,
{\it Spektraldarstellung linearer Transformationen des Hilbertschen Raumes},
Springer, Berlin, 1942.

\item{\hbox to14mm{[SzN47]\hfill}} B. Sz\H{o}kefalvi Nagy,
Perturbations des transformations autoadjointes dans l'espace de Hilbert,
{\it Comment Math. Helv.} {\bf 19} (1946/1947), 347--366.

\item{\hbox to14mm{[Tr]\hfill}} M. Tr\'emon,
Polyn$\hat{\text{\rm o}}$mes de degr\'e minimum connectant deux projections
dans une alg\`ebre de Banach,
{\it Lin. Alg. Appl.} {\bf 64} (1985), 115--132.

\item{\hbox to14mm{[Ze]\hfill}} J. Zem\'anek,
Idempotents in Banach algebras,
{\it Bull. London Math. Soc.} {\bf 11} (1979), 177--183.

\leftskip=0pt

\goodbreak

\enddocument